\documentclass[12pt,a4paper]{article}
\usepackage{graphicx}
\usepackage{subeqnarray}
\begin{document}
%
%
%
%
\def\astrobj#1{#1}
\newenvironment{lefteqnarray}{\arraycolsep=0pt\begin{eqnarray}}
{\end{eqnarray}\protect\aftergroup\ignorespaces}
\newenvironment{lefteqnarray*}{\arraycolsep=0pt\begin{eqnarray*}}
{\end{eqnarray*}\protect\aftergroup\ignorespaces}
\newenvironment{leftsubeqnarray}{\arraycolsep=0pt\begin{subeqnarray}}
{\end{subeqnarray}\protect\aftergroup\ignorespaces}
\newcommand{\diff}{{\rm\,d}}
\newcommand{\pprime}{{\prime\prime}}
\newcommand{\szeta}{\mskip 3mu /\mskip-10mu \zeta}
\newcommand{\FC}{\mskip 0mu {\rm F}\mskip-10mu{\rm C}}
\newcommand{\appleq}{\stackrel{<}{\sim}}
\newcommand{\appgeq}{\stackrel{>}{\sim}}
\newcommand{\Int}{\mathop{\rm Int}\nolimits}
\newcommand{\Nint}{\mathop{\rm Nint}\nolimits}
\newcommand{\range}{{\rm -}}
\newcommand{\displayfrac}[2]{\frac{\displaystyle #1}{\displaystyle #2}}
\def\astrobj#1{#1}
%
\title{Musical intervals under 12-note equal temperament: 
a geometrical interpretation}
\author{
{R.~Caimmi}\footnote{
{\it Physics and Astronomy Department, Padua University, Vicolo Osservatorio
3/2, 35122 Padova, Italy}$\quad$
email: roberto.caimmi@unipd.it~~~
fax: 39-049-8278212}
\footnote{Affiliated up to September 30th 2014. Current status: Studioso
Senior. Current position: in retirement due to age limits.}
, {M$^\circ$ A.~Franzon}\footnote{
{\it Associazione Culturale S.\,Nicol\`o, Piazza Pio X 27,
36043 Camisano Vicentino (VI), Italy}$\quad$
email: albertfranzmuu@gmail.com~~~
}
\footnote{Via Fogazzaro 32, 36043 Camisano Vicentino (VI), Italy}
, {S.~Tognon}\footnote{
{\it FBD House, Bluebell, Dublin 12, Ireland}$\quad$
email: stetogno@gamil.com~~~
fax:+353 (0) 1 478 3574 }
\phantom{agga}
}
%
%
\maketitle
\begin{quotation}
\section*{}
\begin{Large}
\begin{center}

Abstract

\end{center}
\end{Large}
\begin{small}

\noindent\noindent
Musical intervals in multiple of semitones under 12-note equal temperament, or
more specifically
pitch-class subsets of assigned cardinality, $n$, $1\le n\le12$, ($n$-chords)
are conceived as positive integer points,
${\sf P_n}\equiv(\ell_1,\ell_2,...,\ell_n)$,
$\ell_1+\ell_2+...+\ell_n=12$,
within an Euclidean $n$-space,
$\Re^n$.   The number of distinct $n$-chords, $N_{\rm C}(n)$, is inferred
from combinatorics with the extension to $n=0$, involving an Euclidean
0-space, $\Re^0$.
The number of repeating $n$-chords, $\Delta N(n)$, or points which are
turned into themselves
during a circular permutation, $T_n$, of their coordinates, is inferred
from algebraic considerations.  
Finally, the total number of $n$-chords,
$N_{\rm M}(n)=N_{\rm C}(n)+\Delta N(n)$, and the number of $T_n$ set classes,
$\nu_{\rm M}(n)=N_{\rm M}(n)/n$, are determined.   Palindrome and pseudo
palindrome $n$-chords are defined and included among repeating $n$-chords,
with regard to an equivalence relation, $T_n/T_nI$, where reflection is added to
circular permutation.   To this respect, the number of $T_n$ set classes
is inferred concerning palindrome and pseudo palindrome $n$-chords,
$\nu_{\rm P}(n)$, and the remaining $n$-chords,
$\nu_{\rm N}(n)=\nu_{\rm M}(n)-\nu_{\rm P}(n)$, yielding a number of
$T_n/T_nI$ set classes,
$\nu_{\rm Q}(n)=\nu_{\rm N}(n)/2+\nu_{\rm P}(n)$.
The above results are reproduced within the framework of a geometrical
interpretation, where positive integer points related to $n$-chords of
cardinality, $n$, belong to a regular inclined $n$-hedron, $\Psi_{12}^n$, the
vertexes lying on the coordinate axes of a Cartesian orthogonal reference
frame, $({\sf O}x_1x_2...x_n)$, at a distance, $x_i=12$, $1\le i\le n$,
from the origin.   Considering $\Psi_{12}^n$ as special cases of lattice
polytopes,
the number of related nonnegative integer points is also determined for
completeness.   A comparison is performed with the results inferred from group
theory.    The symmetry of the number of $n$-chords, $T_n$ set classes,
$T_n/T_nI$ set classes, with regard to cardinality, is interpreted as
intrinsic to $n$-hedrons and, for this reason, expressed via the binomial
formula.   More generally, the symmetry of the results inferred from the group
theory could be conceived as intrinsic to lattice polytopes in $\Re^n$.

\noindent
{\bf Keywords:} pitch-classes; $n$-chords; $T_n$ set classes; $T_n/T_nI$ set
classes; Euclidean $n$-spaces; $n$-hedrons
\end{small}
\end{quotation}

\section{Introduction} \label{intro}

The question of how many musical intervals in multiples of semitones there are
under 12-note equal
temperament (more specifically, how many pitch-class subsets there are of a
given cardinality, or how many there are with respect to one fixed pitch-class
and further organized under various equivalence relations) is one that has
been answered early and often in the music theory and mathematical literature
e.g., \cite{Loq84}\cite{Bac17}.   The problem has been worked out
independently, in particular the partition into equivalence classes under
transposition (or circular permutation) e.g.,
\cite{Rah80}\cite{Mor87}\cite{Vie93}.   To this respect, it is worth
emphasyzing combinatorial
problems are not essentially musical: the same procedure can be applied e.g.,
for the isomer enumeration in chemistry, for spin analysis in physics, and in
general for the investigation of isomorphism classes of objects e.g.,
\cite{Ker91}\cite{Ker99}.   The most elegant way for solving such problems
is the Polya-Burnside method, which was applied to music theory more than
thirty years ago \cite{Rei85}.

Both pitch-class subsets and $T_n$ set classes of each cardinality from 0
through 12 are
familiar to music theorists but the set-class counts, in absence of a deep
knowledge of the group theory, are performed by use of tables enumerating all
the set classes e.g.,
\cite{For73}\cite{Rah80}\cite{Mor87}\cite{Mor91}\cite{Str05}.
An intermediate use between the two extremes mentioned above
relates to standard techniques in classical combinatorial theory and offers
some simple applications to music theory, including the enumeration of
pitch-class subsets
e.g., \cite{Mes44}\cite{Coh91}\cite{Hoo07}.   In addition, following this line
of thought foreshadows certain aspects of the more difficult work involved in
group theory, and therefore may form a pedagogically benefical bridge to the
advanced material e.g., \cite{Ben06} Chap.\,9 \cite{Hoo07}\cite{Fri99}.

To this respect, the present paper is restricted to the simplest case of $T_n$
and $T_n/T_nI$ set classes, with regard to pitch-class subsets (internal
patterns or internal structures) of cardinality, $n$, where the sum of musical
intervals in multiples of semitones equals 12, or $n$-chords.  Of course,
related results are already known in
the literature e.g., \cite{Rei85}\cite{Ben06} Chap.\,9
\cite{Hoo07}\cite{Fri93}\cite{Fri99}\cite{Rea97}\cite{Jed06}, but the
exposition here is expected
to be more readily accessible to music theory community and, last but not 
least, to interest in group theory by itself.   The current approach is
essentially algebraic and geometric: in short, the paper presents ``an
algorithmic theory,'' one of many possible.   The main steps of the method
may be summarized as follows.

First, $n$-chords, $\{\ell_1,\ell_2,...\ell_n\}$, are
conceived as positive integer points of coordinates,
${\sf P_n}\equiv(\ell_1,\ell_2,...\ell_n)$, with respect to a Cartesian
orthogonal reference frame, $({\sf O}x_1x_2...x_n)$, in an Euclidean
$n$-space, $\Re^n$.   In this view, $n$-chords may be
thought of as made of coordinates.   Then $T_n$ set classes of each
cardinality are partitioned into two main categories, namely set classes where
$n$-chords exhibit distinct e.g., (1,2,3,6) and repeating e.g., (2,4,2,4)
coordinates, respectively.

Second, $n$-chords belonging to the above mentioned categories are
enumerated separately and the amount of related $T_n$ set classes is
determined.

Third, the number of $T_n/T_nI$ set classes of each cardinality is also
determined following a similar procedure.

Fourth, further attention is
devoted to the geometrical interpretation in itself.

The method could, in
principle, be extended to musical intervals in multiples of semitones
under $L$-note (instead of 12-note) equal temperament.

The text is organized as follows.   The first, second and third step outlined
above are developed in different subsections of Section \ref{como}.  The
fourth step is considered in Section \ref{gein}.   The discussion is presented
in Section \ref{disc}.      The conclusion is shown in Section \ref{conc}.
As guidance examples, the method is applied to classical birthday-cake and
necklace problem in Appendix \ref{a:bica} and \ref{a:nela}, respectively.
General properties of $n$-hedrons are outlined in Appendix \ref{a:nhed}.

\section{Enumeration of $n$-chords, $T_n$ and $T_n/T_nI$ set classes}
\label{como}

A pitch-class subset is defined to be a subset of the set of twelve
pitch-classes e.g., \cite{Ben06}
Chap.\,9 \S9.14.   In musical terms, natural numbers within the range,
$1\le n\le12$, could be thought of as
representing musical intervals in multiples of semitones, in the twelve tone
equal tempered octave.   Octave equivalence in the musical scale implies two
notes belong to the same pitch-class if they differ by a whole number of
octaves.
Then addition has an obvious interpretation as addition of musical intervals.
To this respect, an origin must be chosen via one fixed pitch-class.  For
further details, an interested reader is addressed to specific investigations
e.g., \cite{Fri99} or textbooks e.g., \cite{Ben06} Chap.\,9.

Accordingly, pitch-class subsets of cardinality, $n$, are denoted as
$n$-tuples, $\{\ell_1,\ell_2,...,\ell_n\}$, where
$\ell_1$, $\ell_2$, ..., $\ell_n$, are natural numbers.   Let $n$-chords be
defined as pitch-class subsets where the boundary condition:
\begin{equation}
\label{eq:boco}
\ell_1+\ell_2+...+\ell_n=12~~;\qquad1\le n\le12~~;
\end{equation}
is satisfied e.g., \cite{Fri99}.

$T_n$ set classes are obtained by transposition (or circular permutation) as
$\{\ell_1,\ell_2,...,\ell_n\}$, $\{\ell_2,\ell_3,...,\ell_1\}$, ...,
$\{\ell_n,\ell_1,...,\ell_{n-1}\}$.   $T_n/T_nI$ set classes are obtained by
reflection (or order inversion) followed by a transposition.   More
specifically, the application of the pitch-class operator, $T_k$, $k\le n$,
on the $n$-tuple, $\{\ell_1,\ell_2,...,\ell_n\}$, yields
$\{\ell_{k+1},...,\ell_n,\ell_1,...,\ell_k\}$, and the application of the
pitch-class operator, $T_kI$, on the same $n$-tuple, yields a reflection,
$\{\ell_n,\ell_{n-1},...,\ell_1\}$, followed by a transposition,
$\{\ell_k,...,\ell_1,\ell_n,...,\ell_{k+1}\}$.

Let the prime form of a set class be defined as a special $n$-chord within
that
class, for which (i) the last element of the $n$-tuple has the larger value
and, in case of multiplicity, (ii) the first element of the $n$-tuple has the
lower value \cite{Cri79}.   For instance, the prime form of the $T_n/T_nI$ set
class,
\begin{displaymath}
\{1,2,9\},\{2,9,1\},\{9,1,2\},\{1,9,2\},\{9,2,1\},\{2,1,9\};
\end{displaymath}
is $\{1,2,9\}$.   Accordingly, $T_n$ and $T_n/T_nI$ set
classes can be
represented by their prime forms, as cyclic ``adjacent internal array''
(CINT$_1$) \cite{Cri79}.

Let $n$-chords be defined as ``distinct'' and ``repeating'' according if
related
$n$-tuples are different or coinciding, respectively.   Let $T_n$ set classes
be defined as ``distinct'' and ``repeating'' according if related $n$-chords
are distinct or repeating, respectively.   For instance, the $T_n$ set class,
\begin{displaymath}
\{1,2,3,6\},\{2,3,6,1\},\{3,6,1,2\},\{6,1,2,3\};
\end{displaymath}
is made of four distinct 4-chords, while the $T_n$ set class,
\begin{displaymath}
\{1,5,1,5\},\{5,1,5,1\},\{1,5,1,5\},\{5,1,5,1\};
\end{displaymath}
is made of two distinct and two repeating 4-chords.   It is worth
mentioning repeating musical intervals in multiples of semitones have been
studied and used since about 70 years ago \cite{Mes44}.

A method shall be exploited in the following subsections, where distinct and
repeating $n$-chords shall be counted separately to yield the number of
$T_n$ and $T_n/T_nI$ set classes.

\subsection{Enumeration of distinct $n$-chords}
\label{codi}

Aiming to a geometrical interpretation, $n$-tuples representing $n$-chords
shall be considered as coordinates of points,
${\sf P_n}\equiv(\ell_1,\ell_2,...,\ell_n)$,
with respect to a Cartesian orthogonal reference
frame, $({\sf O}\,x_1\,x_2\,...\,x_n)$, in an Euclidean $n$-dimension
hyperspace, or $n$-space, $\Re^n$.   More specifically,
${\sf P_n}$ lies within the positive $2^n$-ant (2-ant is versant, 4-ant is
quadrant, 8-ant is octant, and so on), and
the coordinates are natural numbers linked via Eq.\,(\ref{eq:boco}).

With no loss of generality, the dependent coordinate
may be chosen to be $\ell_n$, as:
\begin{equation}
\label{eq:lod}
\ell_n=12-\ell_1-\ell_2-...-\ell_{n-1}~~;\qquad1\le n\le12~~;
\end{equation}
and the projection of ${\sf P_n}$ onto the principal
$(n-1)$-dimension hyperplane, or $(n-1)$-plane,
$({\sf O}x_1x_2...\,x_{n-1})$, is ${\sf P_{n-1}}\equiv(\ell_1,
\ell_2,...,\ell_{n-1})$.   The knowledge of ${\sf P_{n-1}}$
implies the knowledge of ${\sf P_n}$ via Eq.\,(\ref{eq:lod}).

Given a generic projected point, ${\sf P_{n-1}}\equiv(\ell_1,
\ell_2,...,\ell_{n-2},\ell_{n-1})$, let the conjugate point be defined
as ${\sf Q_{n-1}}\equiv(12-\ell_1, S_{2,n-1},...,S_{n-2,n-1},
\ell_{n-1})$, where, in general, $S_{i,n-k}$ are expressed as:
\begin{leftsubeqnarray}
\slabel{eq:Sija}
&& S_{i,n-k}=\ell_i+\ell_{i+1}+...+\ell_{n-k-1}+\ell_{n-k}~~; \\
\slabel{eq:Sijb}
&& 0<\ell_{n-1}<S_{n-2,n-1}<...<S_{3,n-1}<S_{2,n-1}<12-\ell_1~~; \\
\slabel{eq:Sijc}
&& 0<\ell_i<12~~;\qquad1\le i\le n-1~~;
\label{seq:Sij}
\end{leftsubeqnarray}
and Eq.\,(\ref{eq:Sijb}) follows from (\ref{eq:boco}), (\ref{eq:Sijc}).

\begin{figure*}[t]
\begin{center}      
\includegraphics[scale=0.8]{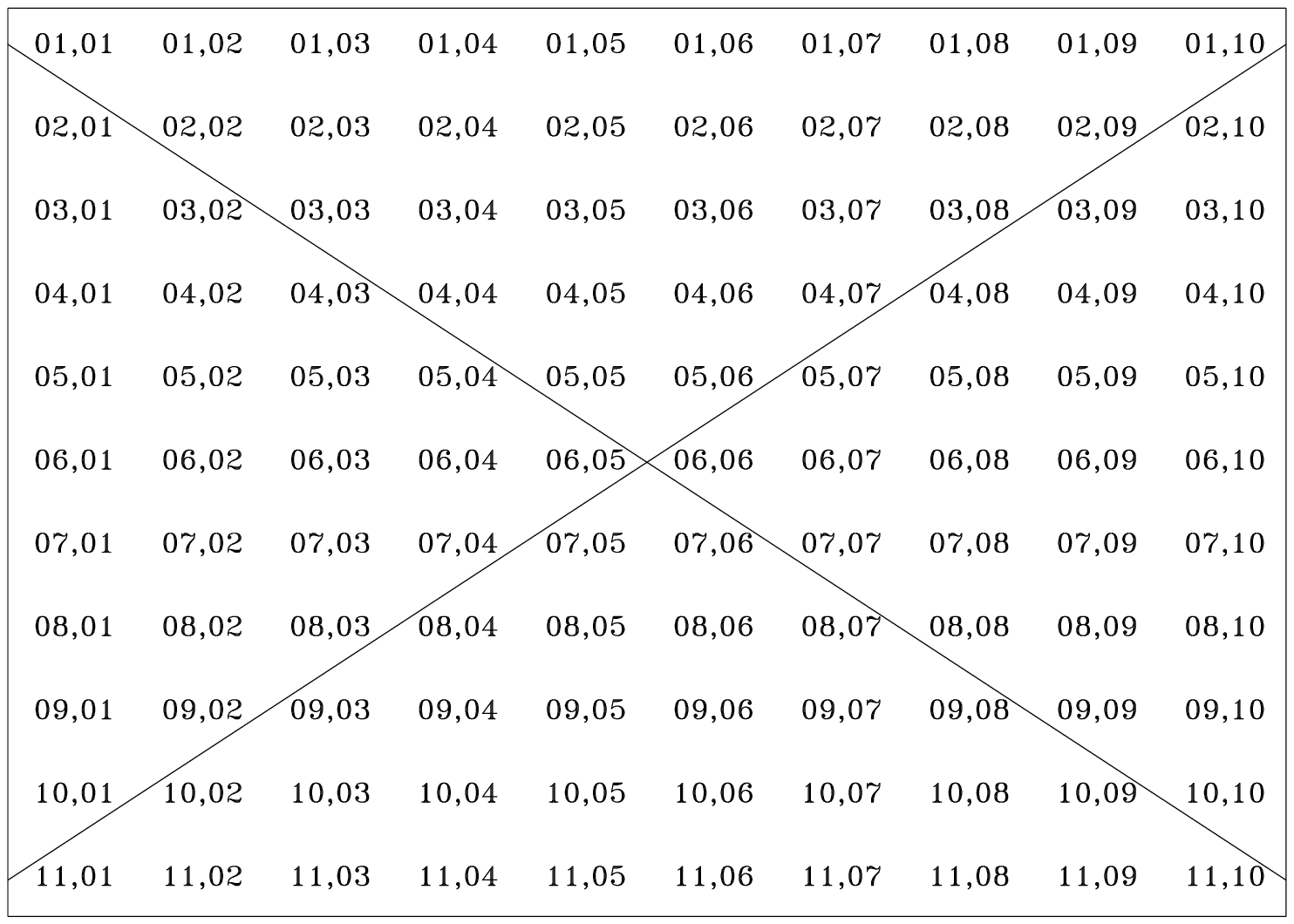}                      
\caption[rbfa]{The $1:1$ correspondence between projected points,
${\sf P_{n-1}}\leftrightarrow{\sf Q_{n-1}}$, in the special case of an
Euclidean $n$-space, $n=3$, ${\sf P_2}(\ell_1,\ell_2)
\leftrightarrow{\sf Q_2}(12-\ell_1,\ell_2)$, yielding a total
number of distinct $n$-chords, $N_{\rm C}(3)=55$.  The
coordinates of projected points, ${\sf P_2}$, satisfying
Eq.\,(\ref{eq:lod}) i.e. $\ell_3=12-\ell_1-\ell_2$, are
placed above the straight line with positive slope.  The
coordinates of projected points, ${\sf Q_2}\leftrightarrow{\sf P_2}$,
are placed below the straight line with negative slope.
The coordinates of projected points, ${\sf P_2}\leftrightarrow
{\sf Q_2}$, ${\sf P_2^\prime}\leftrightarrow{\sf Q_2^\prime}$,
${\sf P_2}\equiv{\sf Q_2^\prime}$, ${\sf P_2^\prime}
\equiv{\sf Q_2}$, are placed within the left angle bisected by a
horizontal straight line.   The coordinates of points for which the
correspondence does not hold and
Eq.\,(\ref{eq:lod}) is violated i.e. $\ell_3>12-\ell_1-\ell_2$,
are placed within the right angle bisected by a horizontal line.
In representing coordinates, brackets have been omitted to save space.
The correspondence, ${\sf P_2}\leftrightarrow{\sf Q_2}$, acts
along columns.   Initial zeroes have been preferred to blank spaces
to save aesthetics.}
\label{f:corn3}
\end{center}       
\end{figure*}

According to Eq.\,(\ref{seq:Sij}), the projected points,
${\sf P_{n-1}}$, ${\sf Q_{n-1}}$, are in a $1:1$ correspondence,
${\sf P_{n-1}}\leftrightarrow{\sf Q_{n-1}}$, or:
\begin{equation}
\label{eq:cPQ}
(\ell_1,\ell_2,...,\ell_{n-2},\ell_{n-1})\leftrightarrow
(12-\ell_1, S_{2,n-1},...,S_{n-2,n-1},\ell_{n-1})~~;
\end{equation}
where the coordinates on the right-hand side of
Eq.\,(\ref{eq:cPQ}) are clearly distinct, monotonically
decreasing, and belonging to the subset of natural numbers,
$\{1,2,...,11\}$, via Eq.\,(\ref{seq:Sij}).   The special
case, $n=3$, is shown in Fig.\,\ref{f:corn3}.
%

With regard to the projected point, ${\sf Q_{n-1}}$, there are
$11-n+n=11-0$ different ways of choosing the first coordinate,
$11-n+(n-1)=11-1$ different ways of choosing the second coordinate
with the preceeding fixed, ..., $11-n+[n-(n-2)]=11-n+2$ different
ways of choosing the $(n-1)$th coordinate with the preceeding fixed,
for a total, $N_{\rm C}^\prime=11\cdot10\cdot...\cdot(11-n+2)=11!/
[11-(n-1)]!$, including points whose coordinates are linked by
permutations i.e. with place exchanged one with respect to the other.

For $(n-1)$ fixed distinct coordinates, there are $(n-1)$ different
ways of choosing the first coordinate, $(n-2)$ different ways of
choosing the second coordinate with the preceeding fixed, ...,
$[n-(n-2)]=2$ different ways of choosing the $(n-2)$th coordinate
with the preceeding fixed, $[n-(n-1)]=1$ univocal way of choosing
the $(n-1)$th coordinate with the preceeding fixed, for a total,
$N^\prime=(n-1)\cdot(n-2)\cdot...\cdot2\cdot1=(n-1)!$.

In conclusion, the total number of projected points, ${\sf Q_{n-1}}$,
having coordinates (i) belonging to the subset of natural numbers,
$\{1,2,...,11\}$; (ii) distinct the one with respect to the other;
(iii) univocally ordered i.e. excluding permutations between coordinates;
is expressed by the ratio, $N_{\rm C}=N_{\rm C}^\prime/N^\prime$, as:
\begin{equation}
\label{eq:NC}
N_{\rm C}=\frac{11!}{[11-(n-1)]!}\frac1{(n-1)!}=\frac n{12}\frac{12!}
{(12-n)!n!}=\frac n{12}{12 \choose n}~~;
\end{equation}
in terms of the binomial coefficients:
\begin{equation}
\label{eq:bico}
{N \choose K}=\frac{N!}{K!(N-K)!}=\frac{N!}{(N-K)!K!}={N \choose N-K}~~;
\end{equation}
related to any pair of natural numbers, $N$, $K$, $N\ge K$.

Accordingly, $N_{\rm C}$ is the number of distinct $n$-chords of
cardinality, $n$.   On the other hand, the total number of pitch-class sets of
cardinality, $n$, regardless of Eq.\,(\ref{eq:boco}), reads $12N_{\rm C}/n$
e.g., \cite{Hoo07} \S29.
Owing to Eq.\,(\ref{eq:bico}), the dependence of
$N_{\rm C}$ on $n$ is symmetric, as shown in Table \ref{t:scale}.
\begin{table}
\caption{Number and fractional number of distinct $n$-chords,
$N_{\rm C}$, $N_{\rm C}/n$, integer part,
$I_{\rm C}=\Int(N_{\rm C}/n)$, number and fractional number of
repeating $n$-chords, $\Delta N$, $\Delta N/n$, total
number and fractional total number of $n$-chords,
$N_{\rm M}=N_{\rm C}+\Delta N$, $N_{\rm M}/n$, for different
cardinality, $n$, $1\le n\le12$.   The additional case, $n=0$, has been added
for completing the symmetry.   See text for further details.}
\label{t:scale}
\begin{center}
\begin{tabular}{|c|c|c|c|c|c|c|c|c|c|c|c|c|c|} \hline
\hline
$n$                        & \phantom{$$}0          & 1 &  2                  &  3                  &   4                  &   5 &   6 &   7 &   8                  &   9                 & 10                  & 11 & 12                     \\
$N_{\rm C}$                & \phantom{$$}0          & 1 & 11                  & 55                  & 165                  & 330 & 462 & 462 & 330                  & 165                 & 55                  & 11 &  1                     \\
$\displayfrac{N_{\rm C}}n$ & $\displayfrac1{12}$    & 1 & $\displayfrac{11}2$ & $\displayfrac{55}3$ & $\displayfrac{165}4$ &  66 &  77 &  66 & $\displayfrac{165}4$ & $\displayfrac{55}3$ & $\displayfrac{11}2$ &  1 & $\displayfrac1{12}$    \\
$I_{\rm C}$                & \phantom{$$}0          & 1 &  5                  & 18                  &  41                  &  66 &  77 &  66 &  41                  &  18                 &  5                  &  1 &  0                     \\
$\Delta N$                 & \phantom{$$}0          & 0 &  1                  &  2                  &  7                   &   0 &  18 &   0 &  14                  &   6                 &  5                  &  0 & 11                     \\
$\displayfrac{\Delta N}n$  & $\displayfrac{11}{12}$ & 0 & $\displayfrac12$    & $\displayfrac23$    & $\displayfrac74$     & 0   &   3 &   0 & $\displayfrac74$     & $\displayfrac23$    & $\displayfrac12$    &  0 & $\displayfrac{11}{12}$ \\
$N_{\rm M}$                & \phantom{$$}0          & 1 & 12                  & 57                  & 172                  & 330 & 480 & 462 & 344                  & 171                 & 60                  & 11 & 12                     \\
$\displayfrac{N_{\rm M}}n$ & \phantom{$$}1          & 1 &  6                  & 19                  &  43                  &  66 &  80 &  66 &  43                  &  19                 &  6                  &  1 &  1                     \\
\hline                            
\end{tabular}                     
\end{center}                      
\end{table}                       
The additional case, $n=0$, has been added for completing the symmetry and
shall be considered below.  Further symmetries
are exhibited by the fractional number of distinct $n$-chords, $N_{\rm C}/n$,
and the related integer part, $I_{\rm C}=\Int(N_{\rm C}/n)$, via
Eq.\,(\ref{eq:bico}).

More specifically, Eq.\,(\ref{eq:NC}) via (\ref{eq:bico}) takes the equivalent
form:
\begin{lefteqnarray}
\label{eq:NCn}
&& \frac{N_{\rm C}}n=\frac1{12}{12 \choose n}=\frac1{12}{12 \choose 12-n}~~;
\end{lefteqnarray}
and the related integer part is:
\begin{equation}
\label{eq:IC}
I_{\rm C}=\Int\left(\frac{N_{\rm C}}n\right)=\Int\left[\frac1{12}
{12\choose n}\right]=\Int\left[\frac1{12}{12 \choose 12-n}\right]~~;
\end{equation}
which is symmetric with respect to the maximum, occurring at $n=6$, as
shown in Table \ref{t:scale}.

To complete the symmetry up to the extreme value, $n=12$, the domain must
be extended down to the opposite extreme, $n=0$, conceived as representing the
empty $n$-chord (no mode).   To this aim, factorials
must be expressed in terms of the Euler Gamma function e.g., \cite{Spi68}
Chap.\,16, as:
\begin{equation}
\label{eq:Gaf}
\Gamma(n+1)=n\Gamma(n)=n!~~;\qquad\Gamma(1)=1~~;\qquad n=1,2,3,...~~;
\end{equation}
where the recursion formula holds for all positive reals, in particular:
\begin{equation}
\label{eq:Ga0}
\lim_{n\to0^+}[n\Gamma(n)]=\lim_{n\to0^+}\Gamma(n+1)=\Gamma(1)=1~~;
\end{equation}
or:
\begin{equation}
\label{eq:Ga02}
\lim_{n\to0^+}\Gamma(n)=\lim_{n\to0^+}\frac1n~~;
\end{equation}
and a similar result is found for $n\to0^-$ extending the recursion
formula, Eq.\,(\ref{eq:Gaf}), to the negative real semiaxis.

In terms of the Euler Gamma function, Eq.\,(\ref{eq:NCn}) reads:
\begin{equation}
\label{eq:NCnG}
\frac{N_{\rm C}}n=\frac1n\frac{\Gamma(12)}{\Gamma(13-n)\Gamma(n)}~~;
\end{equation}
which, for positive infinitesimal $n$, takes the expression:
\begin{equation}
\label{eq:NCG0}
\lim_{n\to0^+}\frac{N_{\rm C}}n=\lim_{n\to0^+}\left[\frac1n\frac{\Gamma(12)}
{\Gamma(13-n)\Gamma(n)}\right]=\frac{\Gamma(12)}{\Gamma(13)}\lim_{n\to0^+}
\left[\frac1n\frac1{\Gamma(n)}\right]~~;
\end{equation}
and the combination of Eqs.\,(\ref{eq:Gaf}), (\ref{eq:Ga0}), (\ref{eq:NCG0}),
yields:
\begin{lefteqnarray}
\label{eq:NCG02}
&& \lim_{n\to0^+}\frac{N_{\rm C}}n=\frac{11!}{12!}=\frac1{12}~~; \\
\label{eq:IC0}
&& \lim_{n\to0^+}\Int\left(\frac{N_{\rm C}}n\right)=0~~;
\end{lefteqnarray}
which completes the symmetry of the fractional number, $N_{\rm C}/n$, and
related integer part, $I_{\rm C}$, with respect to the maximum occurring at
$n=6$.    In authors' opinion, the above considerations add something more to
the bare statement, that $0!=1$ holds by definition.

With regard to a selected Euclidean $n$-space, an integer value of the
fractional number, $N_{\rm C}/n$, makes a necessary (but not sufficient)
condition for a one-to-one correspondence between projected points and
coordinates, ${\sf Q_{n-1}}\leftrightarrow\{s_1,s_2,...,s_{n-1}\}$, where
$s_1=12-\ell_1$; $s_k=S_{k,n-1}$, $2\le k\le n-2$; $s_{n-1}=\ell_{n-1}$.
An inspection of Table \ref{t:scale} shows the necessary condition fails
in several cases, which implies the above mentioned correspondence is not
one-to-one i.e. repeating coordinates, related to repeating $n$-chords,
must also be enumerated.

\subsection{Enumeration of repeating $n$-chords}
\label{coco}

With regard to a selected Euclidean $n$-space and a primitive form,
${\sf P_n}\equiv(\ell_1,\ell_2,...,\ell_n)$, of an assigned $T_n$ set class,
repeating (or transpositionally invariant) $n$-chords,
${\sf P_n^\prime}\equiv(\ell_1^\prime,\ell_2^\prime,...,\ell_n^\prime)$,
${\sf P_n^\pprime}\equiv(\ell_1^\pprime,\ell_2^\pprime,...,\ell_n^\pprime)$,
exhibit identical coordinates,
$\ell_j^\prime=\ell_j^\pprime$, $1\le j\le n$.   More specifically, a
necessary condition for the occurrence of repeating $n$-chords is that
the coordinates, $(\ell_1,\ell_2,...,\ell_n)$, equal one to the other in the
same number i.e. $\ell_{11}=\ell_{21}=...=\ell_{i1}$; $\ell_{12}=\ell_{22}=
...=\ell_{i2}$; ...; $\ell_{1k}=\ell_{2k}=...=\ell_{ik}$; where $1\le ik=n\le
12$.

Accordingly, repeating $n$-chords exhibit identical soloes of
coordinates, $(\ell_1,\ell_1,...,\ell_1)$, or identical duoes of coordinates,
$(\ell_1,\ell_2,\ell_1,\ell_2,...,\ell_1,\ell_2)$, or identical trioes of
coordinates, $(\ell_1,\ell_2,\ell_3,\ell_1,\ell_2,\ell_3,...,\ell_1,\ell_2,
\ell_3)$, or identical quartets of coordinates,
$(\ell_1,\ell_2,\ell_3,\ell_4,
\ell_1,\ell_2,\ell_3,\ell_4,\ell_1,\ell_2,\ell_3,\ell_4)$, or identical
quintets of coordinates, $(\ell_1,\ell_2,\ell_3,\ell_4,\ell_5,\ell_1,\ell_2,
\ell_3,\ell_4,\ell_5)$.   Identical
sextets are not considered in that they yield the chromatic scale,
$\{1,1,1,1,1,1,1,1,1,1,1,1\}$, via Eq.\,(\ref{eq:boco}), and thus reduce to
identical soloes.

For repeating $n$-chords of cardinality, $n$, Eq.\,(\ref{eq:boco}) reduces to:
\begin{equation}
\label{eq:bocr}
\ell_1+\ell_2+...+\ell_i=\frac{12}k=\frac{12}ni~~;
\end{equation}
where $i$ is the number of different coordinates and $k=n/i$ is their
multiplicity.   In any case, the number of repeating $n$-chords, $\Delta N_i(n)$,
to be added to the number of distinct $n$-chords, $N_{\rm C}(n)$, has to be
determined for $T_n$ set classes, while $T_n/T_nI$ set classes shall be
considered afterwards.

Identical soloes of coordinates ($i=1$, $k=n$) cannot occur for $n<2$, and
Eq.\,(\ref{eq:bocr}) reduces to:
\begin{equation}
\label{eq:bors}
\ell_1=\frac{12}n~~;\qquad n\ge2~~;
\end{equation}
which implies the existence of $n$-chords with identical soloes of coordinates
provided the ratio on the right-hand side of Eq.\,(\ref{eq:bors}) is integer.
The related $T_n$ set class is made of $n$ identical singletons of $n$-chords,
one to be counted as
distinct and the remaining $(n-1)$ to be added as repeating.   Accordingly,
the number of repeating $n$-chords reads:
\begin{equation}
\label{eq:DeN1}
\Delta N_1(n)=\zeta(12,n)(n-1)\nu_1(n)~~;\qquad n\ge2~~;
\end{equation}
where $\nu_1(n)$ is the number of $T_n$ set classes including $n$-chords which
satisfy Eq.\,(\ref{eq:bors}), more specifically $\nu_1(n)=1$ for
$n=2,3,4,6,12,$ and $\nu_1(n)=0$ otherwise.
In general the function, $\zeta$, is defined as:
\begin{equation}
\label{eq:zita}
\zeta(m_1,m_2)=\cases{
1~~;\qquad \displayfrac{m_1}{m_2}-\Int\left(\displayfrac{m_1}{m_2}\right)=0~~;
& \cr
1~~;\qquad m_2=0~~; & \cr
0~~;\qquad \displayfrac{m_1}{m_2}-\Int\left(\displayfrac{m_1}{m_2}\right)>0~~;
& \cr
}
\end{equation}
where $\Int(x)$ is the integer part of $x$ and $m_1$, $m_2$, $m_1\ge m_2$, are
natural numbers, $m_1=12$, $m_2=n$, in  the case under discussion.

Identical duoes of coordinates  ($i=2$, $k=n/2$) cannot occur for $n<4$, and
Eq.\,(\ref{eq:bocr}) reduces to:
\begin{equation}
\label{eq:bord}
\ell_1+\ell_2=\frac{24}n~~;\qquad n\ge4~~;
\end{equation}
which implies the existence of $n$-chords with identical duoes of coordinates
provided the ratio on the right-hand side of Eq.\,(\ref{eq:bord}) is integer.
Related $T_n$ set classes are made of $n/2$ identical doublets of $n$-chords,
each one to be counted as
distinct and the others to be added as repeating.   Accordingly,
the number of repeating $n$-chords reads:
\begin{equation}
\label{eq:DeN2}
\Delta N_2(n)=\zeta(24,n)(n-2)\nu_2(n)~~;\qquad n\ge4~~;
\end{equation}
where $\nu_2(n)$ is the number of $T_n$ set classes including $n$-chords which
satisfy Eq.\,(\ref{eq:bord}), to be determined for $n=4,6,8,$ as $\nu_2(n)=0$
otherwise.

For $n=4$, Eq.\,(\ref{eq:bord}) reduces to $\ell_1+\ell_2=6$ which has
distinct (i.e. at least one different from the other) solutions as
$(\ell_1,\ell_2)=(1,5),(2,4),$ implying $\nu_2(4)=1+1=2$,
$\Delta N_2(4)=2+2=4$.

For $n=6$, Eq.\,(\ref{eq:bord}) reduces to $\ell_1+\ell_2=4$ which has
distinct
solutions as $(\ell_1,\ell_2)=(1,3),$ implying
$\nu_2(6)=1$, $\Delta N_2(6)=2+2=4$.

For $n=8$, Eq.\,(\ref{eq:bord}) reduces to $\ell_1+\ell_2=3$ which has
distinct solutions as $(\ell_1,\ell_2)=(1,2),$ implying $\nu_2(8)=1$,
$\Delta N_2(8)=3+3=6$.   The 8-chord, $\{1,2,1,2,1,2,1,2\}$, is quoted among
{\it modes \`a transpositions limit\'ees} \cite{Mes44}.

Identical trioes of coordinates  ($i=3$, $k=n/3$) cannot occur for $n<6$, and
Eq.\,(\ref{eq:bocr}) reduces to:
\begin{equation}
\label{eq:bort}
\ell_1+\ell_2+\ell_3=\frac{36}n~~;\qquad n\ge6~~;
\end{equation}
which implies the existence of $n$-chords with identical trioes of coordinates
provided the ratio on the right-hand side of Eq.\,(\ref{eq:bort}) is integer.
Related $T_n$ set classes are made of $n/3$ identical triplets of $n$-chords,
each one to be counted as
distinct and the others to be added as repeating.   Accordingly,
the number of repeating $n$-chords reads:
\begin{equation}
\label{eq:DeN3}
\Delta N_3(n)=\zeta(36,n)(n-3)\nu_3(n)~~;\qquad n\ge6~~;
\end{equation}
where $\nu_3(n)$ is the number of $T_n$ set classes including $n$-chords which
satisfy Eq.\,(\ref{eq:bort}), to be determined for $n=6,9,$ as $\nu_3(n)=0$
otherwise.

For $n=6$, Eq.\,(\ref{eq:bort}) reduces to $\ell_1+\ell_2+\ell_3=6$ which has
distinct solutions as $(\ell_1,\ell_2,\ell_3)=(1,2,3),(1,1,4),$ implying
$\nu_3(6)=2+1=3$, $\Delta N_3(6)=3+3+3=9$.
The 6-chord, $\{1,4,1,1,4,1\}$, is quoted among {\it modes \`a
transpositions limit\'ees} \cite{Mes44}.

For $n=9$, Eq.\,(\ref{eq:bort}) reduces to $\ell_1+\ell_2+\ell_3=4$ which has
distinct solutions as $(\ell_1,\ell_2,\ell_3)=(1,1,2),$ implying
$\nu_3(9)=1$, $\Delta N_3(9)=2+2+2=6$.
The 9-chord, $\{2,1,1,2,1,1,2,1,1\}$, is quoted among {\it modes \`a
transpositions limit\'ees} \cite{Mes44}.

Identical quartets of coordinates ($i=4$, $k=n/4$) cannot occur for $n<8$, and
Eq.\,(\ref{eq:bocr}) reduces to:
\begin{equation}
\label{eq:borq}
\ell_1+\ell_2+\ell_3+\ell_4=\frac{48}n~~;\qquad n\ge8~~;
\end{equation}
which implies the existence of $n$-chords with identical quartets of
coordinates
provided the ratio on the right-hand side of Eq.\,(\ref{eq:borq}) is integer.
Related $T_n$ set classes are made of $n/4$ identical quadruplets of 
$n$-chords, each one to be counted as
distinct and the others to be added as repeating.   Accordingly,
the number of repeating $n$-chords reads:
\begin{equation}
\label{eq:DeN4}
\Delta N_4(n)=\zeta(48,n)(n-4)\nu_4(n)~~;\qquad n\ge8~~;
\end{equation}
where $\nu_4(n)$ is the number of $T_n$ set classes including $n$-chords which
satisfy Eq.\,(\ref{eq:borq}), to be determined for $n=8,$ as $\nu_4(n)=0$
otherwise.

For $n=8$, Eq.\,(\ref{eq:borq}) reduces to $\ell_1+\ell_2+\ell_3+\ell_4=6$
which has
distinct solutions as $(\ell_1,\ell_2,\ell_3,\ell_4)=(1,1,2,2),(1,1,1,3),$
implying $\nu_4(8)=1+1=2$, $\Delta N_4(8)=4+4=8$.
The 8-chords, $\{2,2,1,1,2,2,1,1\}$, $\{1,1,3,1,1,1,3,1\}$, are quoted among
{\it modes \`a transpositions limit\'ees} \cite{Mes44}.

Identical quintets of coordinates ($i=5$, $k=n/5$) cannot occur for $n<10$,
and Eq.\,(\ref{eq:bocr}) reduces to:
\begin{equation}
\label{eq:born}
\ell_1+\ell_2+\ell_3+\ell_4+\ell_5=\frac{60}n~~;\qquad n\ge10~~;
\end{equation}
which implies the existence of $n$-chords with identical quintets of coordinates
provided the ratio on the right-hand side of Eq.\,(\ref{eq:born}) is integer.
Related $T_n$ set classes are made of $n/5$ identical quintuplets of
$n$-chords, each one to be counted as
distinct and the others to be added as repeating.   Accordingly,
the number of repeating $n$-chords reads:
\begin{equation}
\label{eq:DeN5}
\Delta N_5(n)=\zeta(60,n)(n-5)\nu_5(n)~~;\qquad n\ge10~~;
\end{equation}
where $\nu_5(n)$ is the number of $T_n$ set classes including $n$-chords which
satisfy Eq.\,(\ref{eq:born}), to be determined for $n=10,$ as $\nu_5(n)=0$
otherwise.

For $n=10$, Eq.\,(\ref{eq:born}) reduces to
$\ell_1+\ell_2+\ell_3+\ell_4+\ell_5=6$
which has distinct solutions as
$(\ell_1,\ell_2,\ell_3,\ell_4,\ell_5)=(1,1,1,1,2),$
implying $\nu_5(10)=1$, $\Delta N_5(10)=5$.
The 10-chord, $\{1,1,1,2,1,1,1,1,2,1\}$, is quoted among {\it modes \`a
transpositions limit\'ees} \cite{Mes44}.

In summary, the number of repeating $n$-chords of each cardinality, with
regard to $T_n$ set classes, reads:
\begin{equation}
\label{eq:DeN}
\Delta N(n)=\sum_{i=1}^5\Delta N_i(n)=\sum_{i=1}^5\zeta(12i,n)(n-i)\nu_i(n)~~;
\end{equation}
where the number of repeating $n$-chords including identical soloes $(i=1)$,
$\Delta N_1$, duoes $(i=2)$, $\Delta N_2$, trioes $(i=3)$, $\Delta N_3$,
quartets $(i=4)$, $\Delta N_4$, quintets $(i=5)$, $\Delta N_5$, and the
total, $\Delta N=\sum_i\Delta N_i$, are listed in Table \ref{t:repN}.
\begin{table}
\caption{Number of repeating $n$-chords of each cardinality, $n$, including
identical singletons, $\Delta N_1$; doublets, $\Delta N_2$; triplets,
$\Delta N_3$; quadruplets, $\Delta N_4$; quintuplets, $\Delta N_5$; total
number of repeating $n$-chords, $\Delta N=\sum_i\Delta N_i$; total number of 
distinct $n$-chords, $N_{\rm C}$; and total number of distinct + repeating
$n$-chords,
$N_{\rm M}=N_{\rm C}+\Delta N$.   See text for further details.}
\label{t:repN}
\begin{center}
\begin{tabular}{|c|r|r|r|r|r|r|r|r|r|r|r|r|r|} \hline
\hline
$n$                        & \phantom{$$}0          & 1 &  2                  &  3                  &   4                  &   5 &   6 &   7 &   8                  &   9                 & 10                  & 11 & 12                     \\
$\Delta N_1$               & \phantom{$$}0          & 0 &  1                  &  2                  &   3                  &   0 &   5 &   0 &   0                  &   0                 &  0                  &  0 & 11                     \\
$\Delta N_2$               & \phantom{$$}0          & 0 &  0                  &  0                  &   4                  &   0 &   4 &   0 &   6                  &   0                 &  0                  &  0 &  0                     \\
$\Delta N_3$               & \phantom{$$}0          & 0 &  0                  &  0                  &   0                  &   0 &   9 &   0 &   0                  &   6                 &  0                  &  0 &  0                     \\
$\Delta N_4$               & \phantom{$$}0          & 0 &  0                  &  0                  &   0                  &   0 &   0 &   0 &   8                  &   0                 &  0                  &  0 &  0                     \\
$\Delta N_5$               & \phantom{$$}0          & 0 &  0                  &  0                  &   0                  &   0 &   0 &   0 &   0                  &   0                 &  5                  &  0 &  0                     \\
$\Delta N$                 & \phantom{$$}0          & 0 &  1                  &  2                  &   7                  &   0 &  18 &   0 &  14                  &   6                 &  5                  &  0 & 11                     \\
$N_{\rm C}$                & \phantom{$$}0          & 1 & 11                  & 55                  & 165                  & 330 & 462 & 462 & 330                  & 165                 & 55                  & 11 &  1                     \\
$N_{\rm M}$                & \phantom{$$}0          & 1 & 12                  & 57                  & 172                  & 330 & 480 & 462 & 344                  & 171                 & 60                  & 11 & 12                     \\
\hline                            
\end{tabular}                     
\end{center}                      
\end{table}                       

Similarly, the number of $T_n$ set classes including repeating $n$-chords,
$\nu_i$, the total number of $T_n$ set classes including repeating $n$-chords,
$\nu=\sum_i\nu_i$, the number of $T_n$ set classes including only distinct
$n$-chords, $\nu_{\rm C}$, and the total number of $T_n$ set classes,
$\nu_{\rm M}=N_{\rm M}/n$, are listed in Table \ref{t:repn}.
\begin{table}
\caption{Number of $T_n$ set classes including repeating $n$-chords of each
cardinality, $n$, made of
identical soloes, $\nu_1$; duoes, $\nu_2$; trioes,
$\nu_3$; quartets, $\nu_4$; quintets, $\nu_5$; total number of $T_n$ set
classes including repeating $n$-chords, $\nu=\sum_i\nu_i$; total number of
$T_n$ set classes including only distinct $n$-chords, $\nu_{\rm C}$; and total
number of $T_n$ set classes including distinct + repeating $n$-chords,
$\nu_{\rm M}=\nu_{\rm C}+\nu$.   The $T_0$ set class has been arbitrarily
conceived as including $0$-chords made of (no) identical soloes, to
preserve symmetry in $\nu_{\rm C}$ and $\nu$.   See text for further details.}
\label{t:repn}
\begin{center}
\begin{tabular}{|c|r|r|r|r|r|r|r|r|r|r|r|r|r|} \hline
\hline
$n$                   & \phantom{$$}0          & 1 &  2                  &  3                  &   4                  &   5 &   6 &   7 &   8                  &   9                 & 10                  & 11 & 12                     \\
$\nu_1$               & \phantom{$$}1          & 0 &  1                  &  1                  &   1                  &   0 &   1 &   0 &   0                  &   0                 &  0                  &  0 &  1                     \\
$\nu_2$               & \phantom{$$}0          & 0 &  0                  &  0                  &   2                  &   0 &   1 &   0 &   1                  &   0                 &  0                  &  0 &  0                     \\
$\nu_3$               & \phantom{$$}0          & 0 &  0                  &  0                  &   0                  &   0 &   3 &   0 &   0                  &   1                 &  0                  &  0 &  0                     \\
$\nu_4$               & \phantom{$$}0          & 0 &  0                  &  0                  &   0                  &   0 &   0 &   0 &   2                  &   0                 &  0                  &  0 &  0                     \\
$\nu_5$               & \phantom{$$}0          & 0 &  0                  &  0                  &   0                  &   0 &   0 &   0 &   0                  &   0                 &  1                  &  0 &  0                     \\
$\nu$                 & \phantom{$$}1          & 0 &  1                  &  1                  &   3                  &   0 &   5 &   0 &   3                  &   1                 &  1                  &  0 &  1                     \\
$\nu_{\rm C}$         & \phantom{$$}0          & 1 &  5                  & 18                  &  40                  &  66 &  75 &  66 &  40                  &  18                 &  5                  &  1 &  0                     \\
$\nu_{\rm M}$         & \phantom{$$}1          & 1 &  6                  & 19                  &  43                  &  66 &  80 &  66 &  43                  &  19                 &  6                  &  1 &  1                     \\
\hline                            
\end{tabular}                     
\end{center}                      
\end{table}                       

The number of distinct + repeating $n$-chords of each cardinality,
with respect to $T_n$ set classes, via Eqs.\,(\ref{eq:NC}) and (\ref{eq:DeN})
reads:
\begin{equation}
\label{eq:NMCD}
N_{\rm M}(n)=N_{\rm C}(n)+\Delta N(n)=
\frac n{12}{12 \choose n}+\sum_{i=1}^5\zeta(12i,n)(n-i)\nu_i(n)~~;
\end{equation}
and the number of $T_n$ set classes of each cardinality is:
\begin{equation}
\label{eq:nMCD}
\nu_{\rm M}(n)=\frac{N_{\rm M}(n)}n=\frac{N_{\rm C}(n)+\Delta N(n)}n=
\frac 1{12}{12 \choose n}+\frac1n\sum_{i=1}^5\zeta(12i,n)(n-i)\nu_i(n)~~;
\end{equation}
where an inspection of Table \ref{t:scale} shows that, in general, the number
of $T_n$ set classes including only distinct or repeating $n$-chords is
different from $N_{\rm C}(n)/n$ or $\Delta N(n)/n$, respectively.

The above results complete the calculation of $\Delta N$, $\Delta N/n$,
within the domain, $1\le n\le12$, which allows the knowledge of the total
number of (distinct + repeating) $n$-chords, $N_{\rm M}=N_{\rm C}+\Delta N$,
and the total number of $T_n$ set classes,
$\nu_{\rm M}$, which are also listed in Tables \ref{t:scale}, \ref{t:repN},
\ref{t:repn}.   It can be seen $T_n$ set classes
are symmetric with respect to $n=6$, within the domain,
$1\le n\le11$.   The extension of the domain to $n=0$ can be made demanding
symmetry with respect to $n=12$, which implies the following:
\begin{lefteqnarray}
\label{eq:NMn0}
&& \lim_{n\to0^+}\frac{N_{\rm M}(n)}n=\frac{N_{\rm M}(12)}{12}=1~~; \\
\label{eq:Ddn0}
&& \lim_{n\to0^+}\frac{\Delta N(n)}n=\lim_{n\to0^+}\left(\frac{N_{\rm M}(n)}n-
\frac{N_{\rm C}(n)}n\right)=\frac{N_{\rm M}(12)}{12}-\frac{N_{\rm C}(12)}{12}=
\frac{11}{12}~~;\qquad
\end{lefteqnarray}
as shown in Table \ref{t:scale}.

Following a similar procedure, the number of $T_n/T_nI$ set classes of each
cardinality can also be determined.   In this view, for instance, the
$T_n/T_nI$ set class,
\begin{eqnarray*}
&& \{1,2,3,6\}, \{2,3,6,1\}, \{3,6,1,2\}, \{6,1,2,3\}, \\
&& \{6,3,2,1\}, \{1,6,3,2\}, \{2,1,6,3\}, \{3,2,1,6\};
\end{eqnarray*}
is made of eight distinct 4-chords, while the $T_n/T_nI$ set class,
\begin{eqnarray*}
&& \{1,5,5,1\}, \{5,5,1,1\}, \{5,1,1,5\}, \{1,1,5,5\}, \\
&& \{1,5,5,1\}, \{1,1,5,5\}, \{5,1,1,5\}, \{5,5,1,1\};
\end{eqnarray*}
is made of four distinct and four repeating 4-chords.

With regard to a $n$-chord of cardinality, $n$, let the $n$-chord type be
defined as
$\ell_1^{i_1}\ell_2^{i_2}...\ell_k^{i_k}$, where $i_j$ denotes the
multiplicity of the coordinate, $\ell_j$, $1\le j\le k$, which implies
$i_1+i_2+...+i_k=n$.   Clearly the enumeration of distinct $n$-chords remains
unchanged and the results valid for $T_n$ set classes maintain for $T_n/T_nI$
set classes.   Conversely, the number of repeating $n$-chords is expected to grow
due to larger cardinality, $2n$, of $T_n/T_nI$ set classes, which implies
additional kind of repeating $n$-chords with respect to $T_n$ set classes i.e.
$n$-chords made of identical singletons, duoes, trioes, quartets, quintets of
coordinates.

In this view, let palindrome $n$-chords be defined as invariant with respect to
reflection, and pseudo palindrome $n$-chords as invariant with respect to
reflection after appropriate transposition.   For instance, $\{1,5,5,1\}$ is
palindrome
while $\{5,5,1,1\}$ is pseudo palindrome.   Palindrome and pseudo palindrome
$n$-chords make the additional kinds of repeating $n$-chords, in
connection with
$T_n/T_nI$ set classes.   Then the number of repeating palindrome and pseudo
palindrome $n$-chords of each cardinality, which have not previously
considered, has to be determined.

With regard to $T_n/T_nI$ set classes of each cardinality from 1 to 12 (with
the addition of 0), equivalence classes are made of $2n$ $n$-chords which are
related via circular permutation and reflection.   The total number can be
determined along the following steps.
\begin{trivlist}%
\item[\hspace\labelsep{\bf (i)}]
Start from $T_n$ set classes of
each cardinality.
\item[\hspace\labelsep{\bf (ii)}]
Separate $T_n$ set classes exhibiting neither palindrome nor pseudo palindrome
$n$-chords, let the total number be denoted as $\nu_{\rm N}(n)$, from $T_n$
set classes exhibiting palindrome or pseudo palindrome $n$-chords, let the
total number be denoted as $\nu_{\rm P}(n)$.
\item[\hspace\labelsep{\bf (iii)}]
Determine $\nu_{\rm N}(n)$ and $\nu_{\rm P}(n)$.
\item[\hspace\labelsep{\bf (iv)}]
Calculate the total number of $T_n/T_nI$ set classes as
$\nu_{\rm Q}(n)=\nu_{\rm N}(n)/2+\nu_{\rm P}(n)$, according to the above
considerations.
\end{trivlist}%

Palindrome and pseudo palindrome $n$-chords are necessarily made of pairs of
identical coordinates e.g., $\{1,2,3,3,2,1\}$, for even cardinality, with the
addition of a single coordinate e.g., $\{1,4,2,4,1\}$, for odd cardinality.
Then Eq.\,(\ref{eq:boco}) reduces to:
\begin{leftsubeqnarray}
\slabel{eq:bopa}
&& k_1\ell_1+k_2\ell_2+...+k_i\ell_i=12~~; \\
\slabel{eq:bopb}
&& k_1\ell_1+k_2\ell_2+...+k_i\ell_i+\ell_{i+1}=12~~;
\label{seq:bop}
\end{leftsubeqnarray}
respectively,
where $k_j$, $1\le j\le i$, is the multiplicity of the coordinate, $\ell_j$.
Let $T_n$ set classes made of palindrome and pseudo palindrome $n$-chords be
denoted as $T_{n,\,{\rm P}}$.

For $n=0,1,$ $n$-chords remain unchanged after transposition and/or reflection.
For $n=2$, transposition and reflection are equivalent or, in other words, all
$n$-chords are palindrome or pseudo palindrome.   Accordingly,
$\nu_{\rm Q}(n)=\nu_{\rm M}(n)$, $n<3$, where the number of $T_n$ set classes,
$\nu_{\rm M}(n)$, is listed in Table \ref{t:repn}.

For $n=3$, Eq.\,(\ref{eq:bopb}) reduces to:
\begin{lefteqnarray}
\label{eq:bop3}
&& 2\ell_1+\ell_2=12~~;\qquad n=3~~;
\end{lefteqnarray}
which has solutions as $(\ell_1,\ell_2)=(1,10)$, (2,8), (3,6), (4,4), (5,2),
implying $\nu_{\rm P}(3)=5$,
$\nu_{\rm N}(3)=\nu_{\rm M}(3)-\nu_{\rm P}(3)=19-5=14$;
$\nu_{\rm Q}(3)=\nu_{\rm N}(3)/2+\nu_{\rm P}(3)=14/2+5=7+5=12$.

For $n=4$, Eq.\,(\ref{eq:bopa}) reduces to:
\begin{leftsubeqnarray}
\slabel{eq:bop4a}
&& 2\ell_1+\ell_2+\ell_3=12~~;\qquad n=4~~; \\
\slabel{eq:bop4b}
&& 2\ell_1+2\ell_2=12~~;\qquad n=4~~; \\
\slabel{eq:bop4c}
&& 3\ell_1+\ell_2=12~~;\qquad n=4~~;
\label{seq:bop4}
\end{leftsubeqnarray}
which has solutions as $(\ell_1,\ell_2,\ell_3)=(1,2,8)$, (1,3,7), (1,4,6),
(2,1,7), (2,3,5), (3,1,5), (3,2,4), (4,1,3), (3,3,3); $(\ell_1,\ell_2)=(1,5),$
(2,4); $(\ell_1,\ell_2)=(1,9),$ (2,6); respectively.

Case a yields one kind of $T_{n,\,{\rm P}}$ set classes, namely $\{x,y,x,z\}$.
Accordingly, $\nu_{\rm Pa}(4)=1\cdot9=9$.

Case b yields two kinds of $T_{n,\,{\rm P}}$ set classes, namely
$\{x,x,y,y\}$, $\{x,y,x,y\}$.   Accordingly, $\nu_{\rm Pb}(4)=2\cdot2=4$.

Case c yields one kind of $T_{n,\,{\rm P}}$ set classes, namely $\{x,x,x,y\}$.
Accordingly, $\nu_{\rm Pc}(4)=1\cdot2=2$.

Finally, $\nu_{\rm P}(4)=\sum_i\nu_{\rm Pi}(4)=9+4+2=15$;
$\nu_{\rm N}(4)=\nu_{\rm M}(4)-\nu_{\rm P}(4)=43-15=28$;
$\nu_{\rm Q}(4)=\nu_{\rm N}(4)/2+\nu_{\rm P}(4)=28/2+15=14+15=29$.

For $n=5$, Eq.\,(\ref{eq:bopb}) reduces to:
\begin{leftsubeqnarray}
\slabel{eq:bop5a}
&& 4\ell_1+\ell_2=12~~;\qquad n=5~~; \\
\slabel{eq:bop5b}
&& 2\ell_1+2\ell_2+\ell_3=12~~;\qquad n=5~~; \\
\slabel{eq:bop5c}
&& 3\ell_1+2\ell_2=12~~;\qquad n=5~~;
\label{seq:bop5}
\end{leftsubeqnarray}
which has solutions as $(\ell_1,\ell_2)=(1,8)$, (2,4);
$(\ell_1,\ell_2,\ell_3)=(1,2,6),$ (1,3,4), (1,4,2); $(\ell_1,\ell_2)=(2,3);$
respectively.

Case a yields one kind of $T_{n,\,{\rm P}}$ set classes, namely
$\{x,x,x,x,y\}$.  Accordingly, $\nu_{\rm Pa}(5)=1\cdot2=2$.

Case b yields two kinds of $T_{n,\,{\rm P}}$ set classes, namely
$\{x,y,z,y,x\}$, $\{y,x,z,x,y\}$.   Accordingly, $\nu_{\rm Pb}(5)=2\cdot3=6$.

Case c yields two kinds of $T_{n,\,{\rm P}}$ set classes, namely
$\{x,x,x,y,y\}$, $\{x,y,x,y,x\}$.
Accordingly, $\nu_{\rm Pc}(5)=2\cdot1=2$.

Finally, $\nu_{\rm P}(5)=\sum_i\nu_{\rm Pi}(5)=2+6+2=10$;
$\nu_{\rm N}(5)=\nu_{\rm M}(5)-\nu_{\rm P}(5)=66-10=56$;
$\nu_{\rm Q}(5)=\nu_{\rm N}(5)/2+\nu_{\rm P}(5)=56/2+10=28+10=38$.

For $n=6$, Eq.\,(\ref{eq:bopa}) reduces to:
\begin{leftsubeqnarray}
\slabel{eq:bop6a}
&& 2\ell_1+2\ell_2+2\ell_3=12~~;\qquad n=6~~; \\
\slabel{eq:bop6b}
&& 3\ell_1+3\ell_2=12~~;\qquad n=6~~; \\
\slabel{eq:bop6c}
&& 3\ell_1+2\ell_2+\ell_3=12~~;\qquad n=6~~; \\
\slabel{eq:bop6d}
&& 4\ell_1+2\ell_2=12~~;\qquad n=6~~; \\
\slabel{eq:bop6e}
&& 4\ell_1+\ell_2+\ell_3=12~~;\qquad n=6~~; \\
\slabel{eq:bop6f}
&& 5\ell_1+\ell_2=12~~;\qquad n=6~~;
\label{seq:bop6}
\end{leftsubeqnarray}
which has solutions as $(\ell_1,\ell_2,\ell_3)=(1,2,3)$;
$(\ell_1,\ell_2)=(1,3)$; $(\ell_1,\ell_2,\ell_3)=(1,2,5)$, (2,1,4);
$(\ell_1,\ell_2)=(1,4)$; $(\ell_1,\ell_2,\ell_3)=(1,2,6)$, (1,3,5), (2,1,3);
$(\ell_1,\ell_2)=(1,7)$, $(2,2)$; respectively.

Case a yields six kinds of $T_{n,\,{\rm P}}$ set classes, namely
$\{x,y,x,z,y,z\}$, $\{x,z,x,y,z,y\}$, $\{y,x,y,z,x,z\}$, $\{x,z,y,y,z,x\}$,
$\{y,x,z,z,x,y\}$, $\{z,y,x,x,y,z\}$.
Accordingly, $\nu_{\rm Pa}(6)=6\cdot1=6$.

Case b yields two kinds of $T_{n,\,{\rm P}}$ set classes, namely
$\{x,x,x,y,y,y\}$, $\{x,y,x,y,x,y\}$.   Accordingly,
$\nu_{\rm Pb}(6)=2\cdot1=2$.

Case c yields two kinds of $T_{n,\,{\rm P}}$ set classes, namely
$\{x,x,x,y,z,y\}$, $\{x,z,x,y,x,y\}$.
Accordingly, $\nu_{\rm Pc}(6)=2\cdot2=4$.

Case d yields three kinds of $T_{n,\,{\rm P}}$ set classes, namely
$\{x,x,y,y,x,x\}$, $\{x,y,x,x,y,x\}$, $\{x,x,x,y,x,y\}$.
Accordingly, $\nu_{\rm Pd}(6)=3\cdot1=3$.

Case e yields one kind of $T_{n,\,{\rm P}}$ set classes, namely
$\{x,y,x,x,z,x\}$.   Accordingly,
$\nu_{\rm Pe}(6)=1\cdot3=3$.

Case f yields one kind of $T_{n,\,{\rm P}}$ set classes, namely
$\{x,x,x,x,x,y\}$.
Accordingly, $\nu_{\rm Pf}(6)=1\cdot2=2$.

Finally, $\nu_{\rm P}(6)=\sum_i\nu_{\rm Pi}(6)=6+2+4+3+3+2=20$;
$\nu_{\rm N}(6)=\nu_{\rm M}(6)-\nu_{\rm P}(6)=80-20=60$;
$\nu_{\rm Q}(6)=\nu_{\rm N}(6)/2+\nu_{\rm P}(6)=60/2+20=30+20=50$.

For $n=7$, Eq.\,(\ref{eq:bopb}) reduces to:
\begin{leftsubeqnarray}
\slabel{eq:bop7a}
&& 6\ell_1+\ell_2=12~~;\qquad n=7~~; \\
\slabel{eq:bop7b}
&& 4\ell_1+2\ell_2+\ell_3=12~~;\qquad n=7~~; \\
\slabel{eq:bop7c}
&& 5\ell_1+2\ell_2=12~~;\qquad n=7~~;
\label{seq:bop7}
\end{leftsubeqnarray}
which has solutions as $(\ell_1,\ell_2)=(1,6)$;
$(\ell_1,\ell_2,\ell_3)=(1,2,4),$ (1,3,2); $(\ell_1,\ell_2)=(2,1);$
respectively.

Case a yields one kind of $T_{n,\,{\rm P}}$ set classes, namely
$\{x,x,x,x,x,x,y\}$.  Accordingly, $\nu_{\rm Pa}(7)=1\cdot1=1$.

Case b yields three kinds of $T_{n,\,{\rm P}}$ set classes, namely
$\{x,x,y,z,y,x,x\}$, $\{x,y,x,z,x,y,x\}$, $\{y,x,x,z,x,x,y\}$.   Accordingly,
$\nu_{\rm Pb}(7)=3\cdot2=6$.

Case c yields three kinds of $T_{n,\,{\rm P}}$ set classes, namely
$\{y,x,x,x,x,x,y\}$, $\{x,y,,x,x,y,x\}$, $\{x,x,y,x,y,x,x\}$.
Accordingly, $\nu_{\rm Pc}(7)=3\cdot1=3$.

Finally, $\nu_{\rm P}(7)=\sum_i\nu_{\rm Pi}(7)=1+6+3=10$;
$\nu_{\rm N}(7)=\nu_{\rm M}(7)-\nu_{\rm P}(7)=66-10=56$;
$\nu_{\rm Q}(7)=\nu_{\rm N}(7)/2+\nu_{\rm P}(7)=56/2+10=28+10=38$.

For $n=8$, Eq.\,(\ref{eq:bopa}) reduces to:
\begin{leftsubeqnarray}
\slabel{eq:bop8a}
&& 4\ell_1+4\ell_2=12~~;\qquad n=8~~; \\
\slabel{eq:bop8b}
&& 6\ell_1+\ell_2+\ell_3=12~~;\qquad n=8~~; \\
\slabel{eq:bop8c}
&& 5\ell_1+2\ell_2+\ell_3=12~~;\qquad n=8~~; \\
\slabel{eq:bop8d}
&& 6\ell_1+2\ell_2=12~~;\qquad n=8~~; \\
\slabel{eq:bop8e}
&& 7\ell_1+\ell_2=12~~;\qquad n=8~~;
\label{seq:bop8}
\end{leftsubeqnarray}
which has solutions as $(\ell_1,\ell_2)=(1,2)$;
$(\ell_1,\ell_2,\ell_3)=(1,2,4)$; $(\ell_1,\ell_2,\ell_3)=(1,2,3)$;
$(\ell_1,\ell_2)=(1,3)$; $(\ell_1,\ell_2)=(1,5)$; respectively.

Case a yields six kinds of $T_{n,\,{\rm P}}$ set classes, namely
$\{x,x,x,x,y,y,y,y\}$, $\{x,x,y,y,x,x,y,y\}$, $\{x,y,x,y,x,y,x,y\}$,
$\{x,y,x,y,y,x,y,x\}$, $\{x,x,y,y,x,y,y,x\}$, $\{x,x,y,x,x,y,y,y\}$.
Accordingly, $\nu_{\rm Pa}(8)=6\cdot1=6$.

Case b yields one kind of $T_{n,\,{\rm P}}$ set classes, namely
$\{x,x,x,y,x,x,x,z\}$.   Accordingly, $\nu_{\rm Pb}(8)=1\cdot1=1$.

Case c yields three kinds of $T_{n,\,{\rm P}}$ set classes, namely
$\{x,x,x,y,z,y,x,x\}$, $\{x,x,y,x,z,x,y,x\}$, $\{x,y,x,x,z,x,x,y\}$.
Accordingly, $\nu_{\rm Pc}(8)=3\cdot1=3$.

Case d yields four kinds of $T_{n,\,{\rm P}}$ set classes, namely
$\{x,x,x,x,x,x,y,y\}$, $\{x,x,x,x,x,y,x,y\}$, $\{x,x,x,x,y,x,x,y\}$,
$\{x,x,x,y,x,x,x,y\}$.   Accordingly, $\nu_{\rm Pd}(8)=4\cdot1=4$.

Case e yields one kind of $T_{n,\,{\rm P}}$ set classes, namely
$\{x,x,x,x,x,x,x,y\}$.   Accordingly, $\nu_{\rm Pe}(8)=1\cdot1=1$.

Finally, $\nu_{\rm P}(8)=\sum_i\nu_{\rm Pi}(8)=6+1+3+4+1=15$;
$\nu_{\rm N}(8)=\nu_{\rm M}(8)-\nu_{\rm P}(8)=43-15=28$;
$\nu_{\rm Q}(8)=\nu_{\rm N}(8)/2+\nu_{\rm P}(8)=28/2+15=14+15=29$.

For $n=9$, Eq.\,(\ref{eq:bopb}) reduces to:
\begin{leftsubeqnarray}
\slabel{eq:bop9a}
&& 8\ell_1+\ell_2=12~~;\qquad n=9~~; \\
\slabel{eq:bop9b}
&& 6\ell_1+3\ell_2=12~~;\qquad n=9~~;
\label{seq:bop9}
\end{leftsubeqnarray}
which has solutions as $(\ell_1,\ell_2)=(1,4)$; $(\ell_1,\ell_2)=(1,2);$
respectively.

Case a yields one kind of $T_{n,\,{\rm P}}$ set classes, namely
$\{x,x,x,x,x,x,x,x,y\}$.  Accordingly, $\nu_{\rm Pa}(9)=1\cdot1=1$.

Case b yields four kinds of $T_{n,\,{\rm P}}$ set classes, namely
$\{x,x,x,y,y,y,x,x,x\}$, $\{x,x,y,x,y,x,y,x,x\}$, $\{x,y,x,x,y,x,x,y,x\}$,
$\{y,x,x,x,y,x,x,x,y\}$.   Accordingly, $\nu_{\rm Pb}(9)=4\cdot1=4$.

Finally, $\nu_{\rm P}(9)=\sum_i\nu_{\rm Pi}(9)=1+4=5$;
$\nu_{\rm N}(9)=\nu_{\rm M}(9)-\nu_{\rm P}(9)=19-5=14$;
$\nu_{\rm Q}(9)=\nu_{\rm N}(9)/2+\nu_{\rm P}(9)=14/2+5=7+5=12$.

For $n=10$, Eq.\,(\ref{eq:bopa}) reduces to:
\begin{leftsubeqnarray}
\slabel{eq:bop10a}
&& 8\ell_1+2\ell_2=12~~;\qquad n=10~~; \\
\slabel{eq:bop10b}
&& 9\ell_1+\ell_2=12~~;\qquad n=10~~;
\label{seq:bop10}
\end{leftsubeqnarray}
which has solutions as $(\ell_1,\ell_2)=(1,2)$; $(\ell_1,\ell_2)=(1,4)$;
respectively.

Case a yields five kinds of $T_{n,\,{\rm P}}$ set classes, namely
$\{x,x,x,x,x,x,x,x,y,y\}$, $\{x,x,x,x,x,x,x,y,x,y\}$,
$\{x,x,x,x,x,x,y,x,x,y\}$, $\{x,x,x,x,x,y,x,x,x,y\}$,
$\{x,x,x,x,y,x,x,x,x,y\}$.
Accordingly, $\nu_{\rm Pa}(10)=5\cdot1=5$.

Case b yields one kind of $T_{n,\,{\rm P}}$ set classes, namely
$\{x,x,x,x,x,x,x,x,x,y\}$.   Accordingly, $\nu_{\rm Pb}(10)=1\cdot1=1$.

Finally, $\nu_{\rm P}(10)=\sum_i\nu_{\rm Pi}(10)=5+1=6$;
$\nu_{\rm N}(10)=\nu_{\rm M}(10)-\nu_{\rm P}(10)=6-6=0$;
$\nu_{\rm Q}(10)=\nu_{\rm N}(10)/2+\nu_{\rm P}(10)=0/2+6=0+6=6$.

For $n=11$, Eq.\,(\ref{eq:bopb}) reduces to:
\begin{lefteqnarray}
\label{eq:bop11}
&& 10\ell_1+\ell_2=12~~;\qquad n=11~~;
\end{lefteqnarray}
which has solutions as $(\ell_1,\ell_2)=(1,2)$, yielding 
one kind of $T_{n,\,{\rm P}}$ set classes, namely
$\{x,x,x,x,x,x,x,x,x,x,y\}$.  Accordingly, $\nu_{\rm P}(11)=1\cdot1=1$;
$\nu_{\rm N}(11)=\nu_{\rm M}(11)-\nu_{\rm P}(11)=1-1=0$;
$\nu_{\rm Q}(11)=\nu_{\rm N}(11)/2+\nu_{\rm P}(11)=0/2+1=0+1=1$.

For $n=12$, Eq.\,(\ref{eq:bopa}) reduces to:
\begin{lefteqnarray}
\label{eq:bop12}
&& 12\ell_1=12~~;\qquad n=12~~;
\end{lefteqnarray}
which has solutions as $\ell_1=1$, yielding
one kind of $T_{n,\,{\rm P}}$ set classes, namely
$\{x,x,x,x,x,x,x,x,x,x,x,x\}$.   Accordingly, $\nu_{\rm P}(12)=1\cdot1=1$;
$\nu_{\rm N}(12)=\nu_{\rm M}(12)-\nu_{\rm P}(12)=1-1=0$;
$\nu_{\rm Q}(12)=\nu_{\rm N}(12)/2+\nu_{\rm P}(12)=0/2+1=0+1=1$.

The total number of $T_n$ set classes, $\nu_{\rm M}(n)$ (listed in Table
\ref{t:repn} and repeated for better comparison), palindrome and pseudo
palindrome $T_n$ set classes, $\nu_{\rm P}(n)$, neither palindrome nor pseudo
palindrome $T_n$ set classes, $\nu_{\rm N}(n)=\nu_{\rm M}(n)-\nu_{\rm P}(n)$,
and $T_n/T_nI$ set classes, $\nu_{\rm Q}(n)=\nu_{\rm N}(n)/2+\nu_{\rm P}(n)$,
are listed in Table \ref{t:TTI}.
\begin{table}
\caption{Number of $T_n$ set classes, $\nu_{\rm M}(n)$, palindrome and pseudo
palindrome $T_n$ set classes, $\nu_{\rm P}(n)$, neither palindrome nor pseudo
palindrome $T_n$ set classes, $\nu_{\rm N}(n)=\nu_{\rm M}(n)-\nu_{\rm P}(n)$,
and $T_n/T_nI$ set classes, $\nu_{\rm Q}(n)=\nu_{\rm N}(n)/2+\nu_{\rm P}(n)$,
involving $n$-chords of each cardinality, $n$, $0\le n\le12$.   See text for
further details.}
\label{t:TTI}
\begin{center}
\begin{tabular}{|c|r|r|r|r|r|r|r|r|r|r|r|r|r|} \hline
\hline
$n$                   & \phantom{$$}0          & 1 &  2                  &  3                  &   4                  &   5 &   6 &   7 &   8                  &   9                 & 10                  & 11 & 12                     \\
$\nu_{\rm M}$         & \phantom{$$}1          & 1 &  6                  & 19                  &  43                  &  66 &  80 &  66 &  43                  &  19                 &  6                  &  1 &  1                     \\
$\nu_{\rm P}$         & \phantom{$$}1          & 1 &  6                  &  5                  &  15                  &  10 &  20 &  10 &  15                  &   5                 &  6                  &  1 &  1                     \\
$\nu_{\rm N}$         & \phantom{$$}0          & 0 &  0                  & 14                  &  28                  &  56 &  60 &  56 &  28                  &  14                 &  0                  &  0 &  0                     \\
$\nu_{\rm Q}$         & \phantom{$$}1          & 1 &  6                  & 12                  &  29                  &  38 &  50 &  38 &  29                  &  12                 &  6                  &  1 &  1                     \\
\hline                            
\end{tabular}                     
\end{center}                      
\end{table}                       

An inspection of Table \ref{t:TTI} shows a symmetry with respect to $n=6$.
This is why complementation gives a one to one correspondence between $T_n$
set classes of cardinality, $n$ and $12-n$, respectively, which is preserved
for $T_n/T_nI$ set classes e.g., \cite{Ben06} Chap.\,9 \S9.14.   The same
kind of symmetry is also implicit in the binomial formula, expressed by 
Eq.\,(\ref{eq:bico}). 

\section{Geometrical interpretation}\label{gein}

With regard to an Euclidean $n$-space, $\Re^n$, and a Cartesian orthogonal
reference frame,
$({\sf O}\,x_1\,x_2,...\,x_n)$, the extension of the boundary condition,
expressed by Eq.\,(\ref{eq:boco}), to real numbers, reads:
\begin{equation}
\label{eq:R4}
x_1+x_2+...+x_n=12~~;\quad1\le n\le12~~;
\end{equation}
which represents a $(n-1)$-plane intersecting the coordinate axes at the
points, ${\sf V_k}\equiv(12\delta_{1k},12\delta_{2k},...,12\delta_{nk})$,
$1\le k\le n$, where $\delta_{ik}$ is the Kronecker symbol.

The following properties can be established: (i) the $(n-1)$-plane, expressed
by Eq.\,(\ref{eq:R4}), is normal to the $n$-sector ($n=2,$ bisector; $n=3,$
trisector; and so on) of the positive $2^n$-ant; (ii) the region of
$(n-1)$-plane, bounded by the positive $2^n$-ant, is a regular, inclined
$n$-hedron, $\Psi_{12}^n$, of vertexes, ${\sf V_k}$, $1\le k\le n$; (iii)
special cases are $\Psi_{12}^1$, regular vertex; $\Psi_{12}^2$, regualr side;
$\Psi_{12}^3$, regular triangle; $\Psi_{12}^4$, regular tetrahedron; (iv)  the
orthocentre of $\Psi_{12}^n$, ${\sf H_n}\equiv(12/n,12/n,...,12/n)$, is the
intersection between the $(n-1)$-plane and the $n$-sector of the positive
$2^n$-ant.    For further details, an interested reader is addressed to
Appendix \ref{a:nhed}.

In general, $\Psi_{12}^n$ may be divided into $n$ congruent $n$-hedrons,
$\Psi_{12,i}^n$, $1\le i\le n$, by joining the vertexes with the orthocentre.
More specifically, the orthocentre is the common vertex while the remaining
$(n-1)$ vertexes lie on a $(n-2)$ hyperface, or $(n-2)$-face, of
$\Psi_{12}^n$.

\begin{figure*}[t]
\begin{center}      
\includegraphics[scale=0.8]{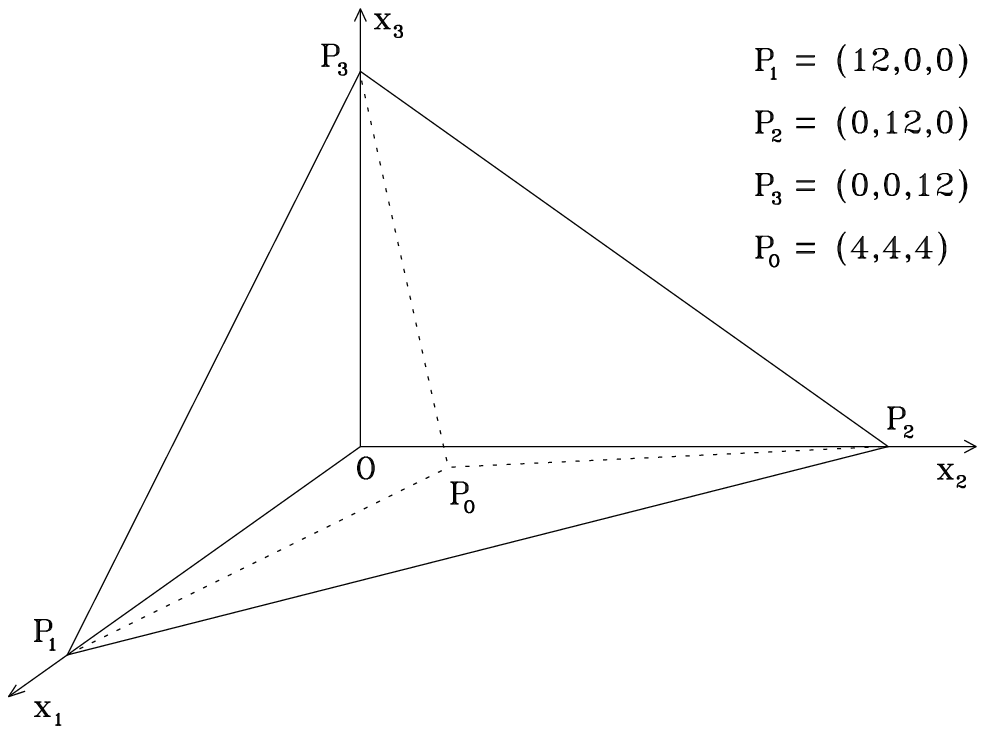}                      
\caption[rbfa]{The special case, $n=3$, where the regular inclined
$n$-hedron reduces to a regular triangle and congruent
$n$-hedrons to isosceles triangles with adjacent equal sides.
The common vertex of isosceles triangles coincides with the orthocentre of the
regular triangle.}
\label{f:ri2t}
\end{center}       
\end{figure*}

\begin{figure*}[t]
\begin{center}      
\includegraphics[scale=0.8]{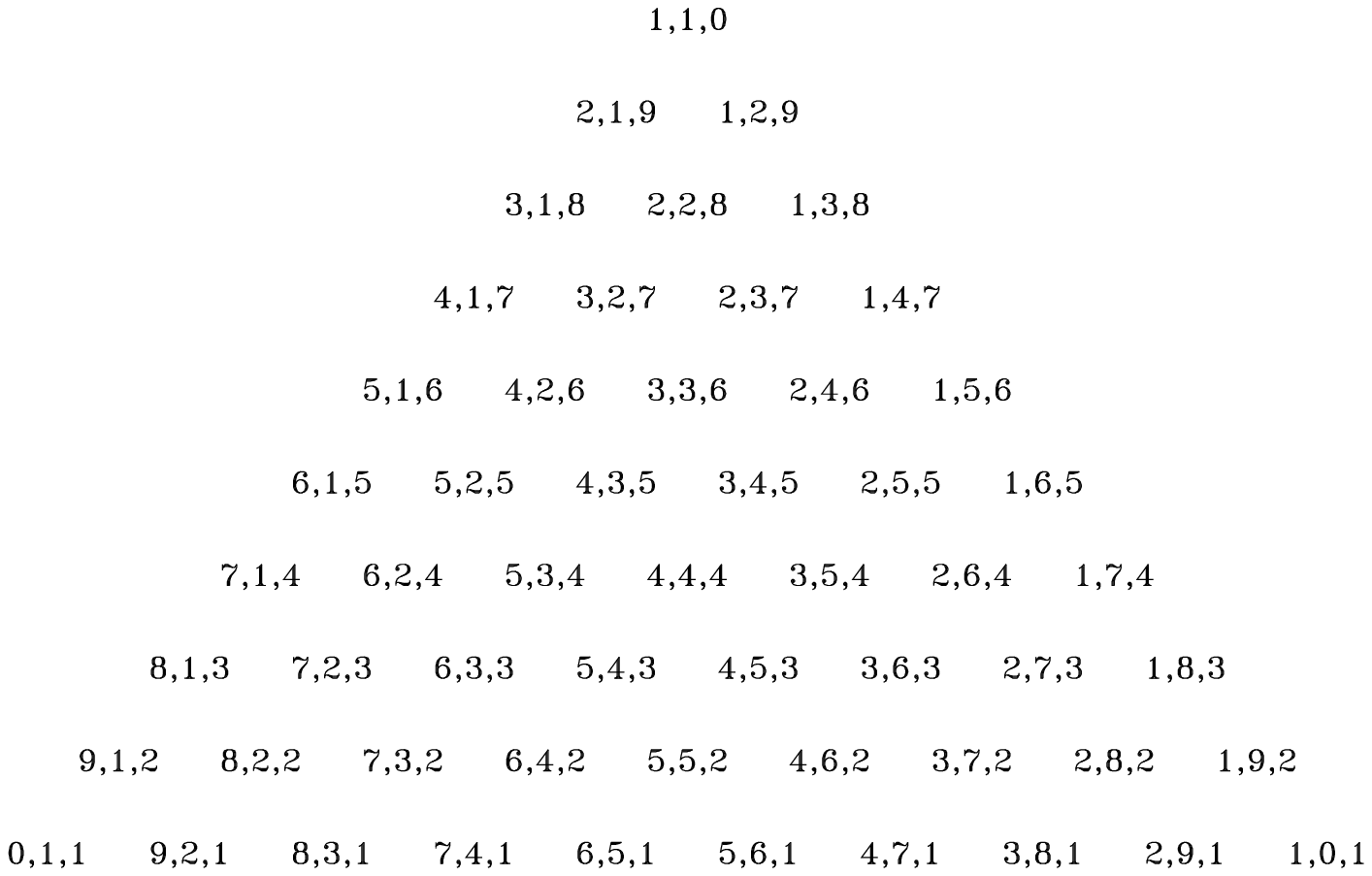}                      
\caption[rbfa]{Representation of distinct $n$-chords of cardinality, $n$,
as positive integer coordinates of points in an
Euclidean $n$-space, lying within a regular inclined $n$-hedron via
Eq.\,(\ref{eq:boco}), in the special case, $n=3$.   The total number of
distinct 3-chords is $N_{\rm C}(3)=55$.  In representing coordinates,
brackets have been omitted to save space and ten has been replaced by zero to
save aesthetics.}
\label{f:rs2t}
\end{center}       
\end{figure*}

The special case,
$\Psi_{12}^3$, is
represented in Fig.\,\ref{f:ri2t} and related $n$-chords are shown as
coordinates of positive integer points, satisfying the boundary condition
expressed by Eq.\,(\ref{eq:boco}), in Fig.\,\ref{f:rs2t}.

According to the above considerations, $n$-chords are represented as
coordinates of positive integer points within $\Psi_{12}^n$ i.e. satisfying
the boundary condition expressed by Eq.\,(\ref{eq:boco}).
$T_n$ and $T_n/T_nI$ set classes contain $n$ and $n+n=2n$ 
points of the kind considered, respectively, which implies (i) points with
distinct coordinates, belonging to the same $T_n$ or $T_n/T_nI$ set class,
are similarly placed within
different $\Psi_{12,i}^n$, $1\le i\le n$; (ii) points with identical soloes,
duoes, trioes, quartets, and quintets of coordinates are placed on
$(n-2)$-faces between different $\Psi_{12,i}^n$, $i=i_1,i_2$,
$1\le i_1<i_2\le n$; (iii) positive integer points (taking into
due account the multiciplity of repeating $n$-chords) are equally partitioned
among $\Psi_{12,i}^n$, $1\le i\le n$, in number of
$\nu_{\rm M}(n)=N_{\rm M}(n)/n$ via Eq.\,(\ref{eq:nMCD}).
For complementary Euclidean $n$-spaces, $\Re^n\leftrightarrow\Re^{12-n}$,
$\nu_{\rm M}(n)=\nu_{\rm M}(12-n)$ as shown in Table \ref{t:repn}.
                                                  
In general, $\Psi_{12,i}^n$, $1\le i\le n$, have basis coinciding with the
$i$th $(n-2)$-face of $\Psi_{12}^n$, which implies $(n-1)$
vertexes in common, with the inclusion of the orthocentre of
$\Psi_{12}^n$.   Accordingly, the coordinates of a generic positive integer
point, ${\sf P_i}\equiv(\ell_{1},\ell_{2},...,\ell_{n})$, belonging to
$\Psi_{12,i}^n$, $1\le i\le n$, have necessarily to satisfy the conditions:
\begin{equation}
\label{eq:cic}
1\le\ell_i\le\frac{12}n~~;\qquad1\le\ell_j\le12-(n-1)~~;\qquad1\le j\le12~~;
\qquad j\ne i~~;
\end{equation}
where the vertex of $\Psi_{12}^n$, placed on the
coordinate axis, $x_i$, does not belong to the congruent $n$-hedron
under consideration.   It is apparent circular permutation of coordinates
makes the points, ${\sf P_i}$, be similarly placed within different
$\Psi_{12,i}^n$, $1\le i\le n$, until the initial configuration is attained.

The extension of the symmetry, outlined in Table \ref{t:scale}, to the special
case, $n=12$, implies Euclidean 0-spaces, $\Re^0$, to be taken into
consideration.
Accordingly, $\Psi_{12}^0$ would lie outside $\Re^0$, on a
coordinate axis at a distance, $x_0=12$, from the origin which, in turn,
coincides with $\Re^0$.   Then the boundary condition, expressed
by Eq.\,(\ref{eq:boco}), is satisfied.

At this stage, a nontrivial question is if $n$-chords of each cardinality can be
enumerated within the framework of the geometrical interpretation outlined
above.   Let $n$-chords exhibiting distinct coordinates be first considered.

For $n=0$, the Euclidean 0-space, $\Re^0$, reduces to the origin and
Eq.\,(\ref{eq:boco}) still holds but with regard to a point of coordinate,
$\ell_0=12$, outside $\Re^0$.   Accordingly, $N_{\rm C}(0)=0$ as listed in
Table \ref{t:scale}.

For $n=1$, the Euclidean 1-space, $\Re^1$, reduces to the real axis, where
Eq.\,(\ref{eq:boco}) is satisfied on the point of coordinate,
$\ell_1=12$.  Accordingly, $N_{\rm C}(1)=1$, as listed in Table \ref{t:scale}.

For $n=2$, the Euclidean 2-space, $\Re^2$, reduces to a plane, where
Eq.\,(\ref{eq:boco}) is satisfied on the points of coordinates,
$(\ell_1,12-\ell_1)$, $1\le\ell_1\le11$.   The above mentioned points are
displaced on a regular inclined 2-hedron (regular side), $\psi_{12}^2$, of
vertexes, ${\sf V_i^\prime}\equiv(1+10\delta_{1i},1+10\delta_{2i})$,
$1\le i\le2$, in number of eleven.
Accordingly, the number of distinct 2-chords must be counted along a
series of superimposed 1-hedrons (regular vertexes) through the second
dimension.  The result is $N_{\rm C}(2)=11$, as listed in Table \ref{t:scale}.

For $n=3$, the Euclidean 3-space, $\Re^3$, is an ordinary space, where
Eq.\,(\ref{eq:boco}) is satisfied on the points of coordinates,
$(\ell_1,\ell_2,12-\ell_1-\ell_2)$, $1\le\ell_k\le10$, $k=1,2$.
As shown in Fig.\,\ref{f:rs2t}, the above mentioned points are
displaced on a regular inclined 3-hedron (regular triangle), $\psi_{12}^3$,
of vertexes,
${\sf V_i^\prime}\equiv(1+9\delta_{1i},1+9\delta_{2i},1+9\delta_{3i})$,
$1\le i\le3$, in number of ten on each 1-face (regular side) and scaled by one
passing to the next related 1-hedron up to the opposite vertex.   Accordingly,
the number of distinct 3-chords must be counted along a series of
superimposed 2-hedrons through the third dimension.   The result is:
\begin{lefteqnarray}
\label{eq:NC3}
&& N_{\rm C}(3)=\sum_{k=1}^{10}k=10+9+8+7+6+5+4+3+2+1=\frac{11\cdot10}2=55\,;
\qquad
\end{lefteqnarray}
as listed in Table \ref{t:scale}.

For $n=4$, Eq.\,(\ref{eq:boco}) is satisfied in $\Re^4$ on the points of
coordinates,
$(\ell_1,\ell_2,\ell_3,$ $12-\ell_1-\ell_2-\ell_3)$, 
$1\le\ell_k\le9$, $k=1,2,3$.   The above mentioned points are
displaced on a regular inclined 4-hedron (regular tetrahedron), $\psi_{12}^4$,
of vertexes,
${\sf V_i^\prime}\equiv(1+8\delta_{1i},1+8\delta_{2i},1+8\delta_{3i},
1+8\delta_{4i})$, $1\le i\le4$, in
number of nine on each 1-face, yielding $10\cdot9/2=45$ points on each 2-face.
Accordingly, the number of distinct 4-chords must be counted along a
series of superimposed 3-hedrons through the fourth dimension.  The result is:
\begin{lefteqnarray}
\label{eq:NC4}
&& N_{\rm C}(4)=\sum_{k=1}^9k+\sum_{k=1}^8k+\sum_{k=1}^7k+\sum_{k=1}^6k+
\sum_{k=1}^5k+\sum_{k=1}^4k+\sum_{k=1}^3k+\sum_{k=1}^2k+\sum_{k=1}^1k
\nonumber \\
&& \phantom{N_{\rm C}(4)}=5\cdot9+4\cdot9+4\cdot7+3\cdot7+3\cdot5+2\cdot5
+2\cdot3+1\cdot3+1\cdot1\qquad \nonumber \\
&& \phantom{N_{\rm C}(4)}=9^2+7^2+5^2+3^2+1^2=165~~;
\end{lefteqnarray}
as listed in Table \ref{t:scale}.

For generic $n$, Eq.\,(\ref{eq:boco}) is satisfied in $\Re^n$ on the points of
coordinates, $(\ell_1,...,\ell_{n-1},12-\ell_1-...-\ell_{n-1})$,
$1\le\ell_k\le12-n+1$, $k=1,2,...,n-1$.    The above mentioned points are
displaced on a regular inclined $n$-hedron,  $\psi_{12}^n$, of vertexes,
${\sf V_i^\prime}\equiv[1+(12-n)\delta_{1i},1+(12-n)\delta_{2i},...,
1+(12-n)\delta_{ni}]$, $1\le i\le n$, in number of $12-n+1$ on each 1-face,
yielding $(12-n+2)(12-n+1)/2$ points on each 2-face.   Accordingly, the number
of distinct $n$-chords must be counted along a series of superimposed
$(n-1)$-hedrons through the $n$th dimension.  The result is:
\begin{lefteqnarray}
\label{eq:NCg}
&& N_{\rm C}(n)=N_{\rm C}\left(\tau_{12-n+1}^{(n-2)}\right)+
N_{\rm C}\left(\tau_{12-n}^{(n-2)}\right)+...
+N_{\rm C}\left(\tau_1^{(n-2)}\right)~~;\qquad n\ge2~~;\qquad
\end{lefteqnarray}
where $N_{\rm C}\left(\tau_k^{(n-2)}\right)$ represents the number of distinct
$n$-chords within the $(n-1)$-hedron of the series, $\tau_k^{(n-2)}$,
$1\le k\le12-n+1$, including on each 1-face
$k$ positive integer points which satisfy Eq.\,(\ref{eq:boco}).   More
specifically, the generic term on the right-hand side of Eq.\,(\ref{eq:NCg})
can be expressed as:
\begin{lefteqnarray}
\label{eq:NCtk}
&& N_{\rm C}\left(\tau_k^{(n-2)}\right)=
N_{\rm C}\left(\tau_{k+1}^{(n-2)}\right)-
N_{\rm C}\left(\tau_{k+1}^{(n-3)}\right)~~;\qquad n\ge2~~;\qquad
\end{lefteqnarray}
where $N_{\rm C}\left(\tau_{k+1}^{(n-3)}\right)$ is the counterpart of
$N_{\rm C}\left(\tau_{k+1}^{(n-2)}\right)$ with regard to the $(n-2)$-hedron
of the related series in $\Re^{n-1}$,
and $N_{\rm C}\left(\tau_{12-n+2}^{(n-2)}\right)=N_{\rm C}(n-1)$.   Then each
term on the right-hand side of Eq.\,(\ref{eq:NCg}) can be determined via
Eq.\,(\ref{eq:NCtk}), provided its counterpart in $\Re^{n-1}$ is known.

The particularization of Eq.\,(\ref{eq:NCtk}) to $n=5$ via
(\ref{eq:NC4}) yields:
\begin{leftsubeqnarray}
\slabel{eq:NC5a}
&& N_{\rm C}\left(\tau_8^{(3)}\right)=N_{\rm C}(4)-
N_{\rm C}\left(\tau_9^{(2)}\right)=165-45=120~~; \\
\slabel{eq:NC5b}
&& N_{\rm C}\left(\tau_7^{(3)}\right)=N_{\rm C}\left(\tau_8^{(3)}\right)-
N_{\rm C}\left(\tau_8^{(2)}\right)=120-36=84~~; \\
\slabel{eq:NC5c}
&& N_{\rm C}\left(\tau_6^{(3)}\right)=N_{\rm C}\left(\tau_7^{(3)}\right)-
N_{\rm C}\left(\tau_7^{(2)}\right)=84-28=56~~; \\
\slabel{eq:NC5d}
&& N_{\rm C}\left(\tau_5^{(3)}\right)=N_{\rm C}\left(\tau_6^{(3)}\right)-
N_{\rm C}\left(\tau_6^{(2)}\right)=56-21=35~~; \\
\slabel{eq:NC5e}
&& N_{\rm C}\left(\tau_4^{(3)}\right)=N_{\rm C}\left(\tau_5^{(3)}\right)-
N_{\rm C}\left(\tau_5^{(2)}\right)=35-15=20~~; \\
\slabel{eq:NC5f}
&& N_{\rm C}\left(\tau_3^{(3)}\right)=N_{\rm C}\left(\tau_4^{(3)}\right)-
N_{\rm C}\left(\tau_4^{(2)}\right)=20-10=10~~; \\
\slabel{eq:NC5g}
&& N_{\rm C}\left(\tau_2^{(3)}\right)=N_{\rm C}\left(\tau_3^{(3)}\right)-
N_{\rm C}\left(\tau_3^{(2)}\right)=10-6=4~~; \\
\slabel{eq:NC5h}
&& N_{\rm C}\left(\tau_1^{(3)}\right)=N_{\rm C}\left(\tau_2^{(3)}\right)-
N_{\rm C}\left(\tau_2^{(2)}\right)=4-3=1~~; 
\label{seq:NC5t}
\end{leftsubeqnarray}
and the particularization of Eq.\,(\ref{eq:NCg}) to $n=5$ via
(\ref{seq:NC5t}) yields the number of distinct 5-chords as:
\begin{lefteqnarray}
\label{eq:NC5}
&& N_{\rm C}(5)=120+84+56+35+20+10+4+1=330~~;
\end{lefteqnarray}
in accordance with Table \ref{t:scale}.

The particularization of Eq.\,(\ref{eq:NCtk}) to $n=6$ via
(\ref{seq:NC5t})-(\ref{eq:NC5}) yields:
\begin{leftsubeqnarray}
\slabel{eq:NC6a}
&& N_{\rm C}\left(\tau_7^{(4)}\right)=N_{\rm C}(5)-
N_{\rm C}\left(\tau_8^{(3)}\right)=330-120=210~~; \\
\slabel{eq:NC6b}
&& N_{\rm C}\left(\tau_6^{(4)}\right)=N_{\rm C}\left(\tau_7^{(4)}\right)-
N_{\rm C}\left(\tau_7^{(3)}\right)=210-84=126~~; \\
\slabel{eq:NC6c}
&& N_{\rm C}\left(\tau_5^{(4)}\right)=N_{\rm C}\left(\tau_6^{(4)}\right)-
N_{\rm C}\left(\tau_6^{(3)}\right)=126-56=70~~; \\
\slabel{eq:NC6d}
&& N_{\rm C}\left(\tau_4^{(4)}\right)=N_{\rm C}\left(\tau_5^{(4)}\right)-
N_{\rm C}\left(\tau_5^{(3)}\right)=70-35=35~~; \\
\slabel{eq:NC6e}
&& N_{\rm C}\left(\tau_3^{(4)}\right)=N_{\rm C}\left(\tau_4^{(4)}\right)-
N_{\rm C}\left(\tau_4^{(3)}\right)=35-20=15~~; \\
\slabel{eq:NC6f}
&& N_{\rm C}\left(\tau_2^{(4)}\right)=N_{\rm C}\left(\tau_3^{(4)}\right)-
N_{\rm C}\left(\tau_3^{(3)}\right)=15-10=5~~; \\
\slabel{eq:NC6g}
&& N_{\rm C}\left(\tau_1^{(4)}\right)=N_{\rm C}\left(\tau_2^{(4)}\right)-
N_{\rm C}\left(\tau_2^{(3)}\right)=5-4=1~~;
\label{seq:NC6t}
\end{leftsubeqnarray}
and the particularization of Eq.\,(\ref{eq:NCg}) to $n=6$ via
(\ref{seq:NC6t}) yields the number of distinct 6-chords as:
\begin{lefteqnarray}
\label{eq:NC6}
&& N_{\rm C}(6)=210+126+70+35+15+5+1=462~~;
\end{lefteqnarray}
in accordance with Table \ref{t:scale}.

The particularization of Eq.\,(\ref{eq:NCtk}) to $n=7$ via
(\ref{seq:NC6t})-(\ref{eq:NC6}) yields:
\begin{leftsubeqnarray}
\slabel{eq:NC7a}
&& N_{\rm C}\left(\tau_6^{(5)}\right)=N_{\rm C}(6)-
N_{\rm C}\left(\tau_7^{(4)}\right)=462-210=252~~; \\
\slabel{eq:NC7b}
&& N_{\rm C}\left(\tau_5^{(5)}\right)=N_{\rm C}\left(\tau_6^{(5)}\right)-
N_{\rm C}\left(\tau_6^{(4)}\right)=252-126=126~~; \\
\slabel{eq:NC7c}
&& N_{\rm C}\left(\tau_4^{(5)}\right)=N_{\rm C}\left(\tau_5^{(5)}\right)-
N_{\rm C}\left(\tau_5^{(4)}\right)=126-70=56~~; \\
\slabel{eq:NC7d}
&& N_{\rm C}\left(\tau_3^{(5)}\right)=N_{\rm C}\left(\tau_4^{(5)}\right)-
N_{\rm C}\left(\tau_4^{(4)}\right)=56-35=21~~; \\
\slabel{eq:NC7e}
&& N_{\rm C}\left(\tau_2^{(5)}\right)=N_{\rm C}\left(\tau_3^{(5)}\right)-
N_{\rm C}\left(\tau_3^{(4)}\right)=21-15=6~~; \\
\slabel{eq:NC7f}
&& N_{\rm C}\left(\tau_1^{(5)}\right)=N_{\rm C}\left(\tau_2^{(5)}\right)-
N_{\rm C}\left(\tau_2^{(4)}\right)=6-5=1~~;
\label{seq:NC7t}
\end{leftsubeqnarray}
and the particularization of Eq.\,(\ref{eq:NCg}) to $n=7$ via
(\ref{seq:NC7t}) yields the number of distinct 7-chords as:
\begin{lefteqnarray}
\label{eq:NC7}
&& N_{\rm C}(7)=252+126+56+21+6+1=462~~;
\end{lefteqnarray}
in accordance with Table \ref{t:scale}.

The particularization of Eq.\,(\ref{eq:NCtk}) to $n=8$ via
(\ref{seq:NC7t})-(\ref{eq:NC7}) yields:
\begin{leftsubeqnarray}
\slabel{eq:NC8a}
&& N_{\rm C}\left(\tau_5^{(6)}\right)=N_{\rm C}(7)-
N_{\rm C}\left(\tau_6^{(5)}\right)=462-252=210~~; \\
\slabel{eq:NC8b}
&& N_{\rm C}\left(\tau_4^{(6)}\right)=N_{\rm C}\left(\tau_5^{(6)}\right)-
N_{\rm C}\left(\tau_5^{(5)}\right)=210-126=84~~; \\
\slabel{eq:NC8c}
&& N_{\rm C}\left(\tau_3^{(6)}\right)=N_{\rm C}\left(\tau_4^{(6)}\right)-
N_{\rm C}\left(\tau_4^{(5)}\right)=84-56=28~~; \\
\slabel{eq:NC8d}
&& N_{\rm C}\left(\tau_2^{(6)}\right)=N_{\rm C}\left(\tau_3^{(6)}\right)-
N_{\rm C}\left(\tau_3^{(5)}\right)=28-21=7~~; \\
\slabel{eq:NC8e}
&& N_{\rm C}\left(\tau_1^{(6)}\right)=N_{\rm C}\left(\tau_2^{(6)}\right)-
N_{\rm C}\left(\tau_2^{(5)}\right)=7-6=1~~;
\label{seq:NC8t}
\end{leftsubeqnarray}
and the particularization of Eq.\,(\ref{eq:NCg}) to $n=8$ via
(\ref{seq:NC8t}) yields the number of distinct 8-chords as:
\begin{lefteqnarray}
\label{eq:NC8}
&& N_{\rm C}(8)=210+84+28+7+1=330~~;
\end{lefteqnarray}
in accordance with Table \ref{t:scale}.

The particularization of Eq.\,(\ref{eq:NCtk}) to $n=9$ via
(\ref{seq:NC8t})-(\ref{eq:NC8}) yields:
\begin{leftsubeqnarray}
\slabel{eq:NC9a}
&& N_{\rm C}\left(\tau_4^{(7)}\right)=N_{\rm C}(8)-
N_{\rm C}\left(\tau_5^{(6)}\right)=330-210=120~~; \\
\slabel{eq:NC9b}
&& N_{\rm C}\left(\tau_3^{(7)}\right)=N_{\rm C}\left(\tau_4^{(7)}\right)-
N_{\rm C}\left(\tau_4^{(6)}\right)=120-84=36~~; \\
\slabel{eq:NC9c}
&& N_{\rm C}\left(\tau_2^{(7)}\right)=N_{\rm C}\left(\tau_3^{(7)}\right)-
N_{\rm C}\left(\tau_3^{(6)}\right)=36-28=8~~; \\
\slabel{eq:NC9d}
&& N_{\rm C}\left(\tau_1^{(7)}\right)=N_{\rm C}\left(\tau_2^{(7)}\right)-
N_{\rm C}\left(\tau_2^{(6)}\right)=8-7=1~~;
\label{seq:NC9t}
\end{leftsubeqnarray}
and the particularization of Eq.\,(\ref{eq:NCg}) to $n=9$ via
(\ref{seq:NC9t}) yields the number of distinct 9-chords as:
\begin{lefteqnarray}
\label{eq:NC9}
&& N_{\rm C}(9)=120+36+8+1=165~~;
\end{lefteqnarray}
in accordance with Table \ref{t:scale}.

The particularization of Eq.\,(\ref{eq:NCtk}) to $n=10$ via
(\ref{seq:NC9t})-(\ref{eq:NC9}) yields:
\begin{leftsubeqnarray}
\slabel{eq:NC10a}
&& N_{\rm C}\left(\tau_3^{(8)}\right)=N_{\rm C}(9)-
N_{\rm C}\left(\tau_4^{(7)}\right)=165-120=45~~; \\
\slabel{eq:NC10b}
&& N_{\rm C}\left(\tau_2^{(8)}\right)=N_{\rm C}\left(\tau_3^{(8)}\right)-
N_{\rm C}\left(\tau_3^{(7)}\right)=45-36=9~~; \\
\slabel{eq:NC10c}
&& N_{\rm C}\left(\tau_1^{(8)}\right)=N_{\rm C}\left(\tau_2^{(8)}\right)-
N_{\rm C}\left(\tau_2^{(7)}\right)=9-8=1~~;
\label{seq:NC10t}
\end{leftsubeqnarray}
and the particularization of Eq.\,(\ref{eq:NCg}) to $n=10$ via
(\ref{seq:NC10t}) yields the number of distinct 10-chords as:
\begin{lefteqnarray}
\label{eq:NC10}
&& N_{\rm C}(10)=45+9+1=55~~;
\end{lefteqnarray}
in accordance with Table \ref{t:scale}.

The particularization of Eq.\,(\ref{eq:NCtk}) to $n=11$ via
(\ref{seq:NC10t})-(\ref{eq:NC10}) yields:
\begin{leftsubeqnarray}
\slabel{eq:NC11a}
&& N_{\rm C}\left(\tau_2^{(9)}\right)=N_{\rm C}(10)-
N_{\rm C}\left(\tau_3^{(8)}\right)=55-45=10~~; \\
\slabel{eq:NC11b}
&& N_{\rm C}\left(\tau_1^{(9)}\right)=N_{\rm C}\left(\tau_2^{(9)}\right)-
N_{\rm C}\left(\tau_2^{(8)}\right)=10-9=1~~;
\label{seq:NC11t}
\end{leftsubeqnarray}
and the particularization of Eq.\,(\ref{eq:NCg}) to $n=11$ via
(\ref{seq:NC11t}) yields the number of distinct 11-chords as:
\begin{lefteqnarray}
\label{eq:NC11}
&& N_{\rm C}(11)=10+1=11~~;
\end{lefteqnarray}
in accordance with Table \ref{t:scale}.

The particularization of Eq.\,(\ref{eq:NCtk}) to $n=12$ via
(\ref{seq:NC11t})-(\ref{eq:NC11}) yields:
\begin{lefteqnarray}
\label{eq:NC12t}
&& N_{\rm C}\left(\tau_1^{(10)}\right)=N_{\rm C}(11)-
N_{\rm C}\left(\tau_2^{(9)}\right)=11-10=1~~;
\end{lefteqnarray}
and the particularization of Eq.\,(\ref{eq:NCg}) to $n=12$ via
(\ref{eq:NC12t}) yields the number of distinct 12-chords as:
\begin{lefteqnarray}
\label{eq:NC12}
&& N_{\rm C}(12)=1~~;
\end{lefteqnarray}
in accordance with Table \ref{t:scale}.

In summary, Eqs.\,(\ref{eq:NC3})-(\ref{eq:NC12}) disclose Eq.\,(\ref{eq:NC})
can be inferred via geometric as well as algebraic considerations, but the
former alternative deserves further attention.   To this respect, let regular,
inclined $n$-hedrons, $\Psi_{12}^n$, of vertexes,
${\sf V_i}\equiv(12\delta_{1i},12\delta_{2i},...,12\delta_{ni})$,
$1\le i\le n$, be considered as special cases of lattice polytopes in $\Re^n$
e.g., \cite{BP99}.

In this view, Eq.\,(\ref{eq:NC}) yields the number of lattice
points (i.e. points of integer coordinates) within $\Psi_{12}^n$, while
the total number reads:
\begin{lefteqnarray}
\label{eq:NV}
&& N_{\rm V}(n)=N_{\rm C}(n)+N_{\rm S}(n)~~;
\end{lefteqnarray}
where $N_{\rm S}(n)$ is the number of lattice points on the $(n-2)$-surface of
$\Psi_{12}^n$, which necessarily exhibit at least one null coordinate.
In the following, $N_{\rm S}(n)$ and $N_{\rm V}(n)$ shall be determined for
each cardinality, keeping in mind values of $N_{\rm C}(n)$ are listed in
Table \ref{t:scale}.

For $n=0$, no lattice polytope exists in the sense $\Psi_{12}^0$ lies outside
$\Re^0$.   Accordingly, $N_{\rm S}(0)=0$, $N_{\rm V}(0)=0$.

For $n=1$, $\Psi_{12}^1$ reduces to a regular vertex,
${\sf V_i}\equiv(12\delta_{1i})$.   Accordingly, $N_{\rm S}(1)=0$,
$N_{\rm V}(1)=N_{\rm C}(1)=1$.

For $n=2$, $\Psi_{12}^2$ reduces to a regular side of extremes,
${\sf V_i}\equiv(12\delta_{1i},12\delta_{2i})$.
Accordingly, $N_{\rm S}(2)=1\cdot2=2$, $N_{\rm V}(2)=11+2=13$.

For $n=3$, $\Psi_{12}^3$ reduces to a regular triangle of vertexes,
${\sf V_i}\equiv(12\delta_{1i},12\delta_{2i},12\delta_{3i})$, as shown in
Fig.\,\ref{f:ri2t}.  Lattice points exhibiting at least one null coordinate
lie on the principal planes, in number of 12 provided vertexes (exhibiting
two null coordinates) are partitioned one per principal plane.   Accordingly,
$N_{\rm S}(2)=3\cdot12=36$, $N_{\rm V}(2)=55+36=91$.

In general, passing from $\Re^{n-1}$ to $\Re^n$, the number of lattice points
within $\Psi_{12}^{n-1}$ has to be incremented by the number of lattice
points within $\Psi_{12}^{n-2}$ for including points exhibiting at least
one null coordinate, which makes the number of lattice points on a
$(n-2)$-face of $\Psi_{12}^n$ i.e. $1/n$ the number of lattice points on
the $(n-2)$-surface of $\Psi_{12}^n$.   Accordingly, the number of
lattice points on the $(n-2)$-surface of $\Psi_{12}^n$ reads:
\begin{lefteqnarray}
\label{eq:NS}
&& N_{\rm S}(n)=n[N_{\rm C}(n-1)+N_{\rm C}(n-2)]~~;\qquad n>1~~;
\end{lefteqnarray}
which, for $n=2,3,$ reproduces the above results.

Using Eqs.\,(\ref{eq:NV})-(\ref{eq:NS}) for each cardinality yields:
\begin{leftsubeqnarray}
\slabel{eq:NSa}
&& N_{\rm S}(2)=2(1+0)=2~~;\qquad N_{\rm V}(2)=11+2=13~~; \\
\slabel{eq:NSb}
&& N_{\rm S}(3)=3(11+1)=36~~;\qquad N_{\rm V}(3)=55+36=91~~; \\
\slabel{eq:NSc}
&& N_{\rm S}(4)=4(55+11)=264~~;\qquad N_{\rm V}(4)=165+264=429~~; \\
\slabel{eq:NSd}
&& N_{\rm S}(5)=5(165+55)=1100~~;\qquad N_{\rm V}(5)=330+1100=1430~~; \\
\slabel{eq:NSe}
&& N_{\rm S}(6)=6(330+165)=2970~~;\qquad N_{\rm V}(6)=462+2970=3432~~; \\
\slabel{eq:NSf}
&& N_{\rm S}(7)=7(462+330)=5544~~;\qquad N_{\rm V}(7)=462+5544=6006~~; \\
\slabel{eq:NSg}
&& N_{\rm S}(8)=8(462+462)=7392~~;\qquad N_{\rm V}(8)=330+7392=7722~~; \\
\slabel{eq:NSh}
&& N_{\rm S}(9)=9(330+462)=7128~~;\qquad N_{\rm V}(9)=165+7128=7293~~; \\
\slabel{eq:NSi}
&& N_{\rm S}(10)=10(165+330)=4950~~;\qquad N_{\rm V}(10)=55+4950=5005~~;~~~ \\
\slabel{eq:NSj}
&& N_{\rm S}(11)=11(55+165)=2420~~;\qquad N_{\rm V}(11)=11+2420=2431~~; \\
\slabel{eq:NSk}
&& N_{\rm S}(12)=12(11+55)=792~~;\qquad N_{\rm V}(12)=1+792=793~~;
\label{seq:NSt}
\end{leftsubeqnarray}
which are listed in Table \ref{t:poli} together with $N_{\rm C}(n)$, the last
repeated for better comparison.
\begin{table}
\caption{Number of lattice points within the $(n-1)$-volume, $N_{\rm C}$, and
on the $(n-2)$-surface, $N_{\rm S}$, of $\Psi_{12}^n$ polytopes, and total
amount, $N_{\rm V}=N_{\rm C}+N_{\rm S}$.   See text for further details.}
\label{t:poli}
\begin{center}
\begin{tabular}{|c|r|r|r|r|r|r|r|r|r|r|r|r|r|} \hline
\hline
$n$                 & \phantom{$$}0          & 1 &  2                  &  3                  &   4                  &    5 &    6 &    7 &    8                  &    9                 &   10                &   11 &  12                     \\
$N_{\rm C}$         & \phantom{$$}0          & 1 & 11                  & 55                  & 165                  &  330 &  462 &  462 &  330                  &  165                 &   55                &   11 &   1                     \\
$N_{\rm S}$         & \phantom{$$}0          & 0 &  2                  & 36                  & 264                  & 1110 & 2970 & 5944 & 7392                  & 7128                 & 4950                & 2420 & 792                     \\
$N_{\rm V}$         & \phantom{$$}0          & 1 & 13                  & 91                  & 429                  & 1430 & 3432 & 6006 & 7722                  & 7293                 & 5005                & 2431 & 793                     \\
\hline                            
\end{tabular}                     
\end{center}                      
\end{table}                       

At this stage, what remains to be done is the enumeration of repeating
$n$-chords
within the framework of the current geometrical interpretation.   To this
respect, it is worth remembering $\Psi_{12}^n$ can be partitioned into $n$
congruent $\Psi_{12,i}^n$, $1\le i\le n$, where different $n$-chords within
an assigned $T_n$ or $T_n/T_nI$ set class relate to points which are similarly
placed within different
$\Psi_{12,i}^n$.   In this view, repeating $n$-chords necessarily lie on the
boundary between different $\Psi_{12,i}^n$ i.e. $(n-2)$-faces, and must be
counted twice, three times, or more, according if points under consideration
are common to two, three, or more $\Psi_{12,i}^n$.
%
%
Then repeating
$n$-chords are expected to relate to some kind of symmetric points, lying on
the straight lines joining the orthocentres of $(j-1)$-faces,
$1\le j\le n-2$.

For $n=0$, $\Psi_{12}^0$ lies outside $\Re^0$, which implies no repeating
0-chords.

For $n=1$, $\Psi_{12}^1=\Psi_{12,1}^1$, which implies no repeating 1-chords.

For $n=2$, $\Psi_{12,1}^2$ and $\Psi_{12,2}^2$ have in common the orthocentre
of $\Psi_{12}^2$, ${\sf H_2}\equiv(6,6)$, which has to be counted twice for
being equally partitioned among $\Psi_{12,1}^2$ and $\Psi_{12,2}^2$, for a
total of repeating 1.   In summary, the number of repeating 2-chords is
1, as shown in Table \ref{t:repN}.

For $n=3$, $\Psi_{12,1}^3$, $\Psi_{12,2}^3$, $\Psi_{12,3}^3$, have in common
the orthocentre of $\Psi_{12}^3$, ${\sf H_3}\equiv(4,4,4)$, which has to be
counted three times for being equally partitioned among $\Psi_{12,1}^3$,
$\Psi_{12,2}^3$, $\Psi_{12,3}^3$, for a total of repeating 2.   In summary,
the number of repeating 3-chords is 2, as shown in Table \ref{t:repN}.

For $n=4$, $\Psi_{12,i}^4$, $1\le i\le4$, have in common
the orthocentre of $\Psi_{12}^4$, ${\sf H_4}\equiv(3,3,3,3)$, which has to be
counted four times for being equally partitioned among $\Psi_{12,i}^4$,
$1\le i\le4$, for a total of repeating 3.    In addition, pairs of 2-faces
(regular triangles) have in common orthocentres of 1-faces (regular sides) on
the boundary, ${\sf H_{4,2}}\equiv[\ell+k(\delta_{1i}+\delta_{1j}),
\ell+k(\delta_{2i}+\delta_{2j}),...,\ell+k(\delta_{4i}+\delta_{4j})]$,
$\ell=1,2$, $k=4,2$.   The number of 1-face orthocentres equals 6 in any case,
where 4 are equally partitioned among $\Psi_{12,i}^4$, while the remaining 2,
$\{x,y,x,y\}$, $\{y,x,y,x\}$, must necessarily be counted twice for a total of
repeating 2+2=4.   In summary, the number of repeating 4-chords is 3+4=7,
as shown in Table \ref{t:repN}.

For $n=5$, the orthocentre of $\Psi_{12}^5$, is not an integer point.   On the
other hand, $(j-1)$-faces of  $\Psi_{12}^5$, $0<j<4$, occur in multiples of 5
as shown in Table \ref{t:face}, which implies related orthocentres are
equally partitioned among $\Psi_{12,i}^5$, $1\le i\le5$, and repeating
5-chords do not take place, according to Table \ref{t:repN}.

For $n=6$, $\Psi_{12,i}^6$, $1\le i\le6$, have in common the orthocentre of
$\Psi_{12}^6$, ${\sf H_6}\equiv(2,2,2,2,2,2)$, which has to be counted six times
for being equally partitioned among $\Psi_{12,i}^6$,
$1\le i\le6$, as shown in Table \ref{t:repN}.    In addition, pairs of 2-faces
(regular triangles) have in common orthocentres of 1-faces (regular sides) on
the boundary, ${\sf H_{6,2}}\equiv[\ell+k(\delta_{1i}+\delta_{1j}),
\ell+k(\delta_{2i}+\delta_{2j}),...,\ell+k(\delta_{6i}+\delta_{6j})]$,
$\ell=1$, $k=3$.   The number of 1-face orthocentres equals 15,
where 12 are equally partitioned among $\Psi_{12,i}^6$, while the remaining 3,
$\{x,x,y,x,x,y\}$, $\{x,y,x,x,y,x\}$, $\{y,x,x,y,x,x\}$, must necessarily be
counted twice
for a total of repeating $1+1+1=3$.   Furthermore, pairs of 3-faces (regular
tetrahedrons) have in common orthocentres of 2-faces (regular triangles) on
the boundary,
${\sf H_{6,3}}\equiv[\ell+k(\delta_{1i}+\delta_{1j}+\delta_{1\ell}),
\ell+k(\delta_{2i}+\delta_{2j}+\delta_{2\ell}),...,
\ell+k(\delta_{6i}+\delta_{6j}+\delta_{6\ell})]$,
$\ell=1$, $k=2$.   The number of 2-face orthocentres equals 20,
where 18 are equally partitioned among $\Psi_{12,i}^6$, while the remaining 2,
$\{x,y,x,y,x,y\}$, $\{y,x,y,x,y,x\}$, must necessarily be counted three
times for a total of repeating $2+2=4$.   Finally, triads of 3-faces (regular
tetrahedrons) have in common symmetric points with respect to the above
mentioned orthocentres, exhibiting identical pairs of coordinates e.g.,
$\{x,x,y,y,z,z\}$, $z=(x+y)/2$.   The total number of pairs within a 6-tuple
is ${6\choose2}=15$, and the total number of pairs within the remaining
4-tuple is ${4\choose2}=6$, while the remaining 2-tuple is fixed.
Accordingly, the
total number of 6-tuples with pairs of identical coordinates equals
$15\cdot6=90$, where 84 are equally partitioned among $\Psi_{12,i}^6$, while
the remaining 6, $\{x,y,z,x,y,z\}$, $\{y,z,x,y,z,x\}$, $\{z,x,y,z,x,y\}$;
$\{x,z,y,x,z,y\}$, $\{z,y,x,z,y,x\}$, $\{y,x,z,y,x,z\}$; must necessarily be
counted twice for a total of repeating $(1+1+1)+(1+1+1)=3+3=6$.   In summary,
the number of repeating $n$-chords is 5+3+4+6=18, as shown in Table
\ref{t:repN}.

For $n=7$, the orthocentre of $\Psi_{12}^7$, is not an integer point.   On the
other hand, $(j-1)$-faces of  $\Psi_{12}^7$, $0<j<6$, occur in multiples of 7
as shown in Table \ref{t:face}, which implies related orthocentres are
equally partitioned among $\Psi_{12,i}^7$, $1\le i\le7$, and repeating
7-chords do not take place, according to Table \ref{t:repN}.

For $n=8$, the orthocentre of $\Psi_{12}^8$, is not an integer point.   On the
other hand, pairs of 2-faces (regular triangles) have in common orthocentres
of 1-faces (regular sides) on
the boundary, ${\sf H_{8,2}}\equiv[\ell+k(\delta_{1i}+\delta_{1j}),
\ell+k(\delta_{2i}+\delta_{2j}),...,\ell+k(\delta_{8i}+\delta_{8j})]$,
$\ell=1$, $k=2$.   The number of 1-face orthocentres equals 28,
where 24 are equally partitioned among $\Psi_{12,i}^8$, while the remaining 4,
$\{x,x,x,y,x,x,x,y\}$, $\{x,x,y,x,x,x,y,x\}$, $\{x,y,x,x,x,y,x,x\}$,
$\{y,x,x,x,y,x,x,x\}$, must necessarily be counted twice for a total of
repeating $1+1+1+1=4$.   In addition, pairs of 4-faces have in common
orthocentres of 3-faces (regular tetrahedrons) on the boundary,
${\sf H_{8,4}}\equiv
[\ell+k(\delta_{1i_1}+\delta_{1i_2}+\delta_{1i_3}+\delta_{1i_4}),
\ell+k(\delta_{2i_1}+\delta_{2i_2}+\delta_{2i_3}+\delta_{2i_4}),...,
\ell+k(\delta_{8i_1}+\delta_{8i_2}+\delta_{8i_3}+\delta_{8i_4})]$,
$\ell=1$, $k=2$.   The number of 3-face orthocentres equals 70,
where 64 are equally partitioned among $\Psi_{12,i}^8$, while the remaining 6,
$\{x,y,x,y,x,y,x,y\}$, $\{y,x,y,x,y,x,y,x\}$; $\{x,x,y,y,x,x,y,y\}$,
$\{x,y,y,x,x,y,y,x\}$, $\{y,y,x,x,y,y,x,x\}$, $\{y,x,x,y,y,x,x,y\}$; must
necessarily be counted four times in the former case and twice in the latter,
for a total of repeating 3+3=6 and 1+1+1+1=4, respectively, yielding 6+4=10.
In summary, the number of
repeating 8-chords is 4+10=14, as shown in Table \ref{t:repN}. 

For $n=9$, the orthocentre of $\Psi_{12}^9$, is not an integer point.   On the
other hand, pairs of 3-faces (regular tetrahedrons) have in common
orthocentres of 2-faces (regular triangles) on the boundary,
${\sf H_{9,3}}\equiv[\ell+k(\delta_{1i}+\delta_{1j}+\delta_{1m}),
\ell+k(\delta_{2i}+\delta_{2j}+\delta_{2m}),...,
\ell+k(\delta_{9i}+\delta_{9j}+\delta_{9m})]$,
$\ell=1$, $k=2$.   The number of 2-face orthocentres equals 84,
where 81 are equally partitioned among $\Psi_{12,i}^9$, while the remaining 3,
$\{x,x,y,x,x,y,x,x,y\}$, $\{x,y,x,x,y,x,x,y,x\}$, $\{y,x,x,y,x,x,y,x,x\}$,
must necessarily be counted three times for a total of repeating 2+2+2=6.
In summary, the number of
repeating 9-chords equals 6, as shown in Table \ref{t:repN}. 

For $n=10$, the orthocentre of $\Psi_{12}^{10}$, is not an integer point.  On
the other hand, pairs of 2-faces (regular triangles) have in common
orthocentres of 1-faces (regular sides) on
the boundary, ${\sf H_{10,2}}\equiv[\ell+k(\delta_{1i}+\delta_{1j}),
\ell+k(\delta_{2i}+\delta_{2j}),...,\ell+k(\delta_{10,i}+\delta_{10,j})]$,
$\ell=1$, $k=1$.   The number of 1-face orthocentres equals 45,
where 40 are equally partitioned among $\Psi_{12,i}^{10}$, while the remaining
5, $\{x,x,x,x,y,x,x,x,x,y\}$, $\{x,x,x,y,x,x,x,x,y,x\}$,
$\{x,x,y,x,x,x,x,y,x,x\}$, $\{x,y,x,x,x,x,y,x,x,x\}$,
$\{y,x,x,x,x,y,x,x,x,x\}$, must necessarily be counted twice for a total of
repeating 1+1+1+1+1=5.   In summary, the number of
repeating 10-chords equals 5, as shown in Table \ref{t:repN}. 

For $n=11$, the orthocentre of $\Psi_{12}^{11}$, is not an integer point.   On
the other hand, $(j-1)$-faces of  $\Psi_{12}^{11}$, $0<j<10$, occur in
multiples
of 11 as shown in Table \ref{t:face}, which implies related orthocentres are
equally partitioned among $\Psi_{12,i}^{11}$, $1\le i\le11$, and repeating
11-chords do not take place, according to Table \ref{t:repN}.

For $n=12$, $\Psi_{12,i}^{12}$, $1\le i\le12$, have in common the orthocentre
of $\Psi_{12}^{12}$, ${\sf H_{12}}\equiv(1,1,...,1)$, which has to be counted
twelve times for being equally partitioned among $\Psi_{12,i}^{12}$,
$1\le i\le12$.  Keeping in mind no additional positive integer point belongs
to $\Psi_{12}^{12}$, the number of repeating 12-chords equals 11, as shown
in Table \ref{t:repN}.

The above results disclose both distinct and repeating $n$-chords and $T_n$
or $T_n/T_nI$ set classes may be enumerated in the light of a geometrical
interpretation, reproducing results inferred via an algebraic approach.   To
this respect, it is worth noticing pairwise $T_n$ and $T_nI$ set classes, made
of repeating but neither palindrome nor pseudo palindrome $n$-chords, are
different, which implies a $2:1$ correspondence between $T_n$ and $T_n/T_nI$
set classes.   Conversely, pairwise $T_n$ and $T_nI$ set classes, made of
repeating palindrome or pseudo palindrome $n$-chords, are coinciding, which
implies a $1:1$ correspondence between $T_n$ and $T_n/T_nI$ set classes.

For further details on $\Psi_{12}^n$ $n$-hedrons, an interested reader is
addressed to Appendix \ref{a:nhed}.

\section{Discussion}\label{disc}

Concerning musical intervals in multiples of semitones under 12-note
equal temperament, counting pitch-class sets
with respect to an assigned equivalence relation lies between two
extreme cases, namely brute-force listing and elegant (but esoteric)
enumeration procedure from group theory.   The method exploited in Section
\ref{como} implies, on one hand, determination of repeating and palindrome or
pseudo palindrome $n$-chords by use, on the other hand, of a less elegant (but
less esoteric) enumeration procedure.   More specifically, $n$-chords of
cardinality, $n$, are related to positive integer points in $\Re^n$, belonging
to a regular inclined $n$-hedron, $\Psi_{12}^n$, where vertexes are placed on
the coordinate axes of a Cartesian orthogonal reference frame,
$({\sf O}x_1\,x_2\,...\,x_n)$, at a distance, $x_i=12$, $1\le i\le 12$, from
the origin.

The method is aimed to help musicians with basic knowledge of algebra,
combinatorics and group theory, but desirous of following a line of thought to
get results.   This is why, say, the binomial theorem has been inferred
instead of
being directly used.   The procedure can be generalized to musical
intervals in multiples of semitones under
$L$-note (instead of 12-note) equal temperament and equivalence relations
other than $T_n$, $T_n/T_nI$, even if expected to be more cumbersome.   Of
course, group theory remains the more powerful tool to this respect, and the
current approach should, ultimately, predispose the reader towards further
deepening on this subject.

According to Eqs.\,(\ref{eq:NC}) and (\ref{eq:DeN}), the number of $T_n$ set
classes of cardinality, $n$, reads:
\begin{lefteqnarray}
\label{eq:nuMc}
&& \nu_{\rm M}(n)=\frac{N_{\rm M}(n)}n=\frac{N_{\rm C}(n)+\Delta N(n)}n
\nonumber \\
&& \phantom{\nu_{\rm M}(n)}=
\frac1{12}{12 \choose n}+\frac1n\sum_{i=1}^{12}\zeta(12i,n)(n-i)\nu_i(n)~~;
\end{lefteqnarray}
where $\zeta(m_1,m_2)$ is defined by Eq.\,(\ref{eq:zita}) and $\nu_i(n)$ is
listed in Table \ref{t:repn}, keeping in mind $i>n$ and/or $n>5$ implies
$\nu_i(n)=0$.

Within the framework of group theory, the result is e.g., \cite{Hoo07} \S72:
\begin{lefteqnarray}
\label{eq:nuMg}
&& \nu_{\rm M}(n)=\frac1{12}\sum_{j=1}^{12}\zeta(12,j)\zeta(n,j)\phi(j)
{12/j \choose n/j}~~;
\end{lefteqnarray}
where $\phi(k)$ is the Euler phi function e.g., \cite{Ben06} Chap.\,9 \S9.7
and the product, $\zeta(12,j)\zeta(n,j)$, by definition equals unity provided
$j$ is a common divisor of $n$ and 12, and equals zero otherwise.   Keeping in
mind $\phi(1)=1$ e.g., \cite{AS72} \S24.3.2, Eq.\,(\ref{eq:nuMg}) may be
cast uder the form:
\begin{lefteqnarray}
\label{eq:nuMh}
&& \nu_{\rm M}(n)=\frac1{12}{12 \choose n}+\frac1{12}\sum_{j=2}^{12}
\zeta(12,j)\zeta(n,j)\phi(j){12/j \choose n/j}~~;
\end{lefteqnarray}
where the first term on the right-hand side of Eq.\,(\ref{eq:nuMh}) equals its
counterpart in Eq.\,(\ref{eq:nuMc}) that is, via Eq.\,(\ref{eq:NCn}), the
fractional number, $N_{\rm C}/n$,  of distinct $n$-chords.   Then the second
term on the right-hand side of Eq.\,(\ref{eq:nuMh}) equals the fractional
number, $\Delta N/n$, of repeating $n$-chords, as listed in Table
\ref{t:core}.
\begin{table}
\caption{Fractional number of repeating $n$-chords, $\Delta N/n$, of each
cardinality, $n$, $1\le n\le12$, inferred from group theory (GT),
Eq.\,(\ref{eq:nuMh}), and from the current paper (CP), Eq.\,(\ref{eq:nuMc}).
For $n=4,6,8$, related addends are different even if the sum yields equal
results.   For $n=0$, Eq.\,(\ref{eq:NCG02}) has also been used.   See
text for further details.}
\label{t:core}
\begin{center}
\begin{tabular}{|r|r|c|c|} \hline
\hline
$n$ & $\phi$ & $\Delta N/n$ (GT) & $\Delta N/n$ (CP) \\
 0  &  1     &         11/12     &       11/12       \\
 1  &  1     &          0/12     &        0/1        \\
 2  &  1     &          6/12     &        1/2        \\
 3  &  2     &          8/12     &        2/3        \\
 4  &  2     &     (15+6)/12     &    (3+4)/4        \\
 5  &  4     &          0/12     &        0/5        \\
 6  &  2     &  (20+12+4)/12     &  (5+4+9)/6        \\
 7  &  6     &          0/12     &        0/7        \\
 8  &  4     &     (15+6)/12     &    (6+8)/8        \\
 9  &  6     &          8/12     &        6/9        \\
10  &  4     &          6/12     &        5/10       \\
11  & 10     &          0/12     &        0/11       \\
12  &  4     &         11/12     &       11/12       \\
\hline                            
\end{tabular}                     
\end{center}                      
\end{table}                       

Leaving aside the special case of null cardinality, $n=0$, an inspection of
Table \ref{t:core} shows the terms of the sum on the right-hand side of
Eqs.\,(\ref{eq:nuMh}) and (\ref{eq:nuMc}), respectively, coincide with the
exception of $n=4,6,8$, even if an equal result is attained, which implies
related formulations are not intrinsically equivalent.

For $n=0$, $\lim_{n\to0^+}\Delta N(n)/n=11/12$ is inferred from Table
\ref{t:scale} via Eqs.\,(\ref{eq:NCG02}) and (\ref{eq:nuMc}).   On the other
hand, (1+2+2+2+4)/12=11/12 via Eq.\,(\ref{eq:nuMh}), in agreement with the
results listed in Table \ref{t:core}.

A similar comparison could be performed with regard to $T_n/T_nI$ set classes
by use of the counterparts of Eqs.\,(\ref{eq:nuMc}) and (\ref{eq:nuMg}) e.g.,
\cite{Hoo07} \S72 which implies a more complex formulation and for this reason
it has not been considered.

Symmetries found in enumerating $n$-chords, as shown by results listed in
Tables
\ref{t:scale}, \ref{t:repN}, \ref{t:repn}, are intrinsic to the geometrical
interpretation.   More specifically, $n$-hedrons exhibit $(j-1)$-faces and
$(n-j-1)$-faces, $1\le j\le n-1$, in equal number, as listed in Table
\ref{t:face}, which can be expressed via binomial coefficients.   In addition,
the boundary condition expressed by Eq.\,(\ref{eq:boco}) implies regular,
inclined $n$-hedrons, $\Psi_{12,i}^n$ ,where vertexes lie on the coordinate
axes at a distance, $x_i=12$, $1\le i\le n$, from the origin.

It is worth emphasizing symmetries are intrinsic to $n$-hedrons and, for this
reason, can be described using the binomial formula.   In particular, regular
inclined $n$-hedrons, $\Psi_{12}^n$, can be partitioned into $n$ congruent
$n$-hedrons, $\Psi_{12,i}^n$, $1\le i\le n$, by joining the orthocentre with
each vertex.   Different $\Psi_{12,i}^n$ have in common $(n-2)$-faces where a
vertex is the orthocentre of $\Psi_{12}^n$.

All of $n$-chords belonging to an assigned $T_n$ or $T_n/T_nI$ set class
are similarly or similarly and symmetrically placed, respectively, within
different $\Psi_{12,i}^n$, where distinct $n$-chords relate to positive
integer points
belonging to $\Psi_{12,i}^n$, $1\le i\le n$, while repeating $n$-chords,
if any, are placed on $(n-2)$-faces of $\Psi_{12,i}^n$ where a vertex is the
orthocentre of $\Psi_{12}^n$.

The geometric representation exploited in the current paper looks
appealing, in that the idea is very simple and very effective, which is
always a good point. The fact that the dimension of the $(n-1)$-plane where
$n$-chords are located only depends on the scale size, and not on the
chromatic gamut in which the scale is embedded (it only affects the number of
points considered) makes it very useful.   One could for example easily
visualise triads in quarter-tone scales, or on even finer-grained scales,
on an ordinary (2-) plane.

Last but not least, the above mentioned geometric
representation has musical significance: rotations around the $n$-sector (i.e.
an axis normal to the $(n-1)$-plane through the orthocentre of $\Psi_{12}^n$)
correspond to chord inversions (in musical
language), or cyclic permutations (in mathematical language). Looking at their
location tells if a scale may have internal symmetries (repeating chords).
Similarly, rotations about an axis joining orthocentres of opposite $(j-1)$
and $(n-1-j)$-faces correspond to chord reflections, which are symmetrically
placed with respect to the rotation axis.
But surely there are many more which would be worth investigating further.

Enumeration of $n$-chords within an algebraic or geometrical framework yields
coinciding results but, in the latter alternative, further perspectives can be
exploited such as counting integer points within lattice polytopes.

Following a similar line of thought, the results of the current paper could be
extended to musical intervals in multiples of semitones under
$L$-note equal temperament, reproducing
results found within the framework of group theory e.g., \cite{Fri99}
\cite{Hoo07}.

The
group theory, of course, remains the more powerful tool in dealing with any
kind of equivalence relation, but further insight on specific problems could
be gained in the light of both algebraic and geometric method outlined in the
current paper.

More specifically, using group theory in counting specified pitch-class sets
could be
related to the description of a high-population $(N\gg1)$ statistical system
e.g., a perfect gas within a box, in terms of observables e.g., volume,
pressure, temperature.   Conversely, the enumeration of pitch-class sets
within an algebraic or geometric framework, as outlined in the current paper,
could be related to the description of a low-population $(N\appgeq1)$
statistical system e.g., perfectly elastic dissipationless spheres moving on a
horizontal plane bounded by perfectly elastic dissipationless walls, in terms
of integration of the motion equations from the knowledge of initial
conditions.

\section{Conclusion}\label{conc}

The current investigation is restricted to musical intervals in
multiples of semitones under 12-note equal temperament, in particular about
the question of how many $n$-chords there are of
a given cardinality, and how many $T_n$ or $T_n/T_nI$ set classes.   In the
light of a geometrical interpretation, $n$-chords of a given cardinality,
$\{\ell_1,\ell_2,...,\ell_n\}$, satisfying the boundary condition expressed by
Eq.\,(\ref{eq:boco}), are related to positive integer points in $\Re^n$,
belonging to a regular inclined $n$-hedron, $\Psi_{12}^n$.   Joining vertexes
with the orthocentre yields $n$ congruent $n$-hedrons, $\Psi_{12,i}^n$,
$1\le i\le n$, where points with
coordinates related by circular permutations i.e. the $n$ elements of $T_n$
set classes, are equally partitioned and similarly placed.   If $T_n$
set classes exhibit repeating $n$-chords, then related points are lying on
$(n-2)$-faces of $\Psi_{12,i}^n$ where a vertex is the orthocentre of
$\Psi_{12}^n$.

The number, $N_{\rm C}(n)$, of distinct $n$-chords is determined within an 
algebraic framework, and a symmetry is found around $q=13/2$ as
listed in Table \ref{t:scale}.   The fractional
number, $N_{\rm C}(n)/n$, of distinct $n$-chords shows a symmetry around
$n=6$, as
listed in Table \ref{t:scale},
provided the domain is extended to $n=0$ in connection with Euclidean
0-spaces, $\Re^0$.   The fractional numbers, $\Delta N(n)/n$ and
$N_{\rm M}(n)/n=[N_{\rm C}(n)+\Delta N(n)]/n$, of repeating and total
(distinct + repeating) $n$-chords, respectively, also exhibit a symmetry
around $n=6$, including $n=0$, as
listed in Table \ref{t:scale}.
It is worth remembering $\nu_{\rm M}(n)=N_{\rm M}(n)/n$ is the number of
$T_n$ set classes.

Following the same line of thought, the above results can be generalized to
musical intervals in multiples of semitones under $L$-note equal
temperament.   In particular,
Eqs.\,(\ref{eq:NC}), (\ref{eq:NCn}), (\ref{eq:IC}), (\ref{eq:NCnG}),
(\ref{eq:NCG0}), (\ref{eq:NCG02}), hold provided $12\mp k$, $k=0,1$,
is replaced by $L\mp k$ therein.

The above mentioned symmetries are disclosed within a geometrical framework,
where $n$-chords together with $T_n$ and $T_n/T_nI$ set classes can be
enumerated as well.   More
specifically, symmetries are intrinsic to $n$-hedrons, in general on one
hand, and via the boundary condition of the problem on the other hand.

The main results of the current paper can be summarized as follows.
\begin{trivlist}%
\item[\hspace\labelsep{\bf (1)}]
Using a one-to-one correspondence between $n$-chords and
positive integer points within an Euclidean $n$-space, $\Re^n$, both distinct
and repeating $n$-chords are enumerated within the framework of an algebraic
method together with $T_n$ and $T_n/T_nI$ set classes.   Related results equal
their counterparts already known in the context of group theory e.g.,
\cite{Fri99} \cite{Ben06} Chap.\,9 \cite{Hoo07}.
\item[\hspace\labelsep{\bf (2)}]
Positive integer points related to $n$-chords of cardinality, $n$, belong to a
regular inclined $n$-hedron, $\Psi_{12}^n$, of vertexes,
${\sf V_i}\equiv(12\delta_{1i}, 12\delta_{2i},...,12\delta_{ni})$,
$1\le i\le n$.  By joining the orthocentre,  
${\sf H_n}\equiv(12/n, 12/n,...,12/n)$, with vertexes, ${\sf V_i}$,
$\Psi_{12}^n$ can be subdivided into $n$ congruent $n$-hedrons,
$\Psi_{12,i}^n$, $1\le i\le n$.   Positive integer points related to 
$n$-chords within an assigned $T_n$ or $T_n/T_nI$ set class are equally
partitioned among and similarly or similarly and symmetrically placed within
$\Psi_{12,i}^n$.   In particular, points related to repeating $n$-chords are
placed on $(n-2)$-faces
between different $\Psi_{12,i}^n$, to be counted according to their
multiplicity i.e. the number of involved $\Psi_{12,i}^n$.
\item[\hspace\labelsep{\bf (3)}]
Both distinct and repeating $n$-chords are enumerated within the framework of a
geometrical method together with $T_n$ and $T_n/T_nI$ set classes, reproducing
their counterparts inferred from an algebraic method in (1).
\item[\hspace\labelsep{\bf (4)}]
Results in (3) allow the calculation of nonnegative integer (i.e. exhibiting
nonnegative integer coordinates) points belonging to $\Psi_{12}^n$ which, to
this respect, can be considered as special cases of lattice polytopes in
$\Re^n$.
\item[\hspace\labelsep{\bf (5)}]
Symmetries shown by the number of $n$-chords and $T_n$ or $T_n/T_nI$ set classes
with respect
to the cardinality, $n=6$, are intrinsic to the geometry in a twofold manner.
In general, the number of $(j-1)$-faces of $n$-hedrons, $1\le j\le n-1$, can
be expressed via binomial coefficients as listed in Table \ref{t:face}.   In
particular, the boundary condition of the problem, expressed by
Eq.\,(\ref{eq:boco}), implies consideration of regular inclined $n$-hedrons,
$\Psi_{12}^n$, where the number of positive integer points can also be
expressed via the binomial formula, as shown in Table \ref{t:scale}.   To this
respect, $\Re^n$ and $\Re^{12-n}$, $0\le n\le12$, or $\Re^{6-m}$, $\Re^{6+m}$,
$0\le m\le6$, can be considered as Euclidean complementary spaces.
More generally, the symmetry of the results inferred from the group
theory could be conceived as intrinsic to lattice polytopes in $\Re^n$.
\end{trivlist}%

\section*{Acknowledgements}
The authors are indebted to two anonymous referees of an earlier version of
the paper (2008), for critical comments and useful suggestions.

\appendix
\section*{Appendix}

\section{Birthday-cake problem}
\label{a:bica}

Let $N$ candles of $L$ different colour be put at equal distance from the
symmetry axis of a circular birthday cake, on the $N$ vertexes of a regular
polygon.
Accordingly, a generic configuration is defined by a $N$-tuple of colours,
selected among the available $L$, in brief a colour set.   Then $T_n$ set
classes are defined by $N$ repeated rotations of the cake about its axis by an
angle, $\alpha_{\rm N}=2\pi/N$.

The problem is: how many colour sets and
$T_n$ set classes do exist?   For sake of clarity, let attention
be restricted to a simple case, $N=4$ and $L=3$: blue (B), green (G), red (R),
say e.g., \cite{Hoo07} \S46.

A solution to the problem shall be exploited along the following steps.
\begin{trivlist}%
\item[\hspace\labelsep{\bf (i)}]
Calculate the total number of distinct colour sets, 4-tuples in the case under
consideration.
\item[\hspace\labelsep{\bf (ii)}]
Select different kinds of colour sets with regard to the elements of the 
4-tuples.
\item[\hspace\labelsep{\bf (iii)}]
For each kind of colour sets infer the number of repeating 4-tuples within
related $T_n$  set classes.
\item[\hspace\labelsep{\bf (iv)}]
Calculate the total number of distinct + repeating colour sets.
\item[\hspace\labelsep{\bf (v)}]
Calculate the total number of $T_n$ set classes.
\end{trivlist}%

Different configurations of the birthday cake can be related to different
colour sets.   For each place of a 4-tuple, there are three possibilities for
the choice of a colour (B,G,R), yielding a total of $3^4=81$ distinct
4-tuples.   In addition, the following types of colour sets can be defined:
$x^4$, where all colours are identical; $x^3y$, where three colours are
identical and the remaining one is different; $x^2y^2$, where two different
colours
are identical in pairs; $x^2yz$, where two colours are identical and the
remaining two are different.   In general, $x={\rm B,G,R}$; $x\ne y\ne z$.

For $x^4$ colour sets, $\{x,x,x,x\}$, the number of distinct 4-tuples is 
$4!/(4\cdot3\cdot2)=1$, and the number of different elements, $x$, equals
${3\choose1}=3$.   Accordingly, the total number of distinct and repeating
4-tuples is $3\cdot1=3$ and $3\cdot3=9$, respectively.   Then the total number
of distinct + repeating colour sets is $3+9 = 12$.

For $x^3y$ colour sets, $\{x,x,x,y\}$, the number of distinct 4-tuples is 
$4!/(3\cdot2)=4$, and the number of different pairs, $xy$, equals
${4\choose2}=6$.   Accordingly, the total number of distinct and repeating
4-tuples is $6\cdot4=24$ and $0\cdot4=0$, respectively.   Then the total
number of distinct + repeating colour sets is $24+0=24$.

For $x^2y^2$ colour sets, $\{x,x,y,y\}$, the number of distinct 4-tuples is 
$4!/(2\cdot2)=6$, and the number of different pairs, $xy$, equals
${3\choose2}=3$.   Accordingly, the total number of distinct and repeating
4-tuples is $3\cdot6=18$ and $3\cdot(1+1)=6$, respectively, the last due to
$\{x,y,x,y\}$ and $\{y,x,y,x\}$.   Then the total number
of distinct + repeating colour sets is $18+6 = 24$.

For $x^2yz$ colour sets, $\{x,x,y,z\}$, the number of distinct 4-tuples is 
$4!/2=12$, and the number of different triads, $xyz$, equals
${3\choose1}=3$.   Accordingly, the total number of distinct and repeating
4-tuples is $3\cdot12=36$ and $0\cdot12=0$, respectively.   Then the total
number of distinct + repeating colour sets is $36+0=36$.

More specifically, in the case under discussion the total number of colour
sets belonging to an assigned type can be determined along the following
steps.

\begin{trivlist}%
\item[\hspace\labelsep{\bf (a)}]
Start from an assigned $T_n$ set class  of the type considered.
\item[\hspace\labelsep{\bf (b)}]
Exchange $y$ with $z$ (if both present) to get a distinct (if any) $T_n$ set
class.
\item[\hspace\labelsep{\bf (c)}]
Exchange $x$ with $y$ (only once in case of multiplicity) to get a distinct
(if any) $T_n$ set class.
\end{trivlist}%

For $x^4$ colour sets, the application of the method yields:
\begin{lefteqnarray*}
&& \{x,x,x,x\}, \{x,x,x,x\}, \{x,x,x,x\}, \{x,x,x,x\}~~;
\end{lefteqnarray*}
where no additional exchange has to be performed as both $y$ and $z$ are
absent, yielding 1 $T_n$ set class made of 1 distinct and 3 repeating colour
sets.   Keeping in mind $x={\rm B,G,R}$, there are $3\cdot1=3$ $T_n$ set
classes for a total of $1\cdot3=3$ distinct and $3\cdot3=9$ repeating colour
sets, according to the above results.

For $x^3y$ colour sets, the application of the method yields:
\begin{lefteqnarray*}
&& \{x,x,x,y\}, \{x,x,y,x\}, \{x,y,x,x\}, \{y,x,x,x\}~~;
\end{lefteqnarray*}
where no additional exchange has to be performed as $z$ is
absent, yielding 1 $T_n$ set class made of 4 distinct colour
sets.   Keeping in mind $xy={\rm BG,GR,RB,GB,RG,BR}$, ($xy\ne yx$ due to 
different multiplicities in $x$ and $y$), there are $6\cdot1=6$ $T_n$ set
classes for a total of $6\cdot4=24$ distinct and $0\cdot4=0$ repeating colour
sets, according to the above results.

For $x^2y^2$ colour sets, the application of the method yields:
\begin{lefteqnarray*}
&& \{x,y,y,x\}, \{y,y,x,x\}, \{y,x,x,y\}, \{x,x,y,y\}~~; \\
&& \{x,y,x,y\}, \{y,x,y,x\}, \{x,y,x,y\}, \{y,x,y,x\}~~;
\end{lefteqnarray*}
where no additional exchange has to be performed as $z$ is absent, yielding 2
$T_n$ set classes made of (4+2) distinct and 2 repeating colour sets.
Keeping in mind $xy={\rm BG,GR,RB}$, ($xy=yx$ do to equal multiplicities
in $x$ and $y$), there are $3\cdot2=6$ $T_n$ set
classes for a total of $3\cdot(4+2)=18$ distinct and $3\cdot2=6$ repeating
colour sets, according to the above results.

For $x^2yz$ colour sets, the application of the method yields:
\begin{lefteqnarray*}
&& \{x,x,y,z\}, \{x,y,z,x\}, \{y,z,x,x\}, \{z,x,x,y\}~~; \\
&& \{x,x,z,y\}, \{x,z,y,x\}, \{z,y,x,x\}, \{y,x,x,z\}~~; \\
&& \{x,y,x,z\}, \{y,x,z,x\}, \{x,z,x,y\}, \{z,x,y,x\}~~;
\end{lefteqnarray*}
where no additional exchange has to be performed, yielding 3
$T_n$ set classes made of $3\cdot4=12$ distinct colour
sets.   Keeping in mind $x={\rm B,G,R}$; $x\ne y\ne z$;
there are $3\cdot3=9$ $T_n$ set
classes for a total of $3\cdot12=36$ distinct and $0\cdot12=0$ repeating
colour sets, according to the above results.

Finally, the number of distinct, $N_{\rm C}$, repeating, $\Delta N$,
total, $N_{\rm M}$, colour sets, respectively, reads:
\begin{lefteqnarray}
\label{eq:NCb}
&& N_{\rm C}=N_{\rm C}(x^4)+N_{\rm C}(x^3y)+N_{\rm C}(x^2y^2)+N_{\rm C}(x^2yz)
=3+24+18+36 \nonumber \\
&& \phantom{N_{\rm C}=N_{\rm C}(x^4)+N_{\rm C}(x^3y)+N_{\rm C}(x^2y^2)+
N_{\rm C}(x^2yz)}=81~~; \\
\label{eq:DNb}
&& \Delta N=\Delta N(x^4)+\Delta N(x^3y)+\Delta N(x^2y^2)+\Delta N(x^2yz)=
9+0+6+0 \nonumber \\
&& \phantom{\Delta N=\Delta N(x^4)+\Delta N(x^3y)+\Delta N(x^2y^2)+\Delta N
(x^2yz)}=15~~;\qquad \\
\label{eq:NMb}
&& N_{\rm M}=N_{\rm C}+\Delta N=12+24+24+36=96~~; 
\end{lefteqnarray}
and the number of $T_n$ set classes, $\nu_{\rm M}$, is:
\begin{lefteqnarray}
\label{eq:nMb}
&& \nu_{\rm M}=\frac{N_{\rm M}}4=\frac{12}4+\frac{24}4+\frac{24}4+\frac{36}4=
3+6+6+9=24~~;
\end{lefteqnarray}
or 96/4=24.

Within the framework of group theory, the application of Burnside's lemma to
birthday-cake problem yields e.g., \cite{Hoo07} \S46:
\begin{lefteqnarray}
\label{eq:nMB}
&& \nu_{\rm M}=\frac{81}4+\frac34+\frac94+\frac34=\frac{96}4=24~~;
\end{lefteqnarray}
which is in agreement with Eq.\,(\ref{eq:nMb}) but the addends are different.

Within the framework of group theory, the application of Polya's theorem to
birthday-cake problem yields e.g., \cite{Hoo07} \S62:
\begin{lefteqnarray}
\label{eq:nMP}
&& \nu_{\rm M}=(1+1+1)+(1+1+1+1+1+1)+2(1+1+1)+3(1+1+1) \nonumber \\
&& \phantom{\nu_{\rm M}}=3+6+6+9=24~~;
\end{lefteqnarray}
which is in agreement with Eq.\,(\ref{eq:nMb}), exhibiting equal addends.

The result from Eq.\,(\ref{eq:nMB}) is based on Burnside's lemma and does not
involve weights on $T_n$ set classes.   On the other hand, the weighted method
given by Polya's theorem is simply a refinement of the same idea, showing
exactly how the 24 $T_n$ set classes are distributed among the various weights
e.g., \cite{Hoo07} \S63.   For further details, an interested reader is
addressed to specialized monographs e.g., \cite{Hoo07} or textbooks e.g.,
\cite{Ben06} Chap.\,9.   The current method exactly reproduces the results
inferred via Polya's theorem.

\section{Necklace problem}
\label{a:nela}

Let $N$ beads of $L$ different colour be put at equal distance from the
symmetry axis of a circular necklace, on the $N$ vertexes of a regular
polygon.
Accordingly, a generic configuration is defined by a $N$-tuple of colours,
selected among the available $L$, in brief a colour set.   Then $T_n$ set
classes are defined by $N$ repeated rotations of the necklace about its axis
by an angle, $\alpha_{\rm N}=2\pi/N$; $T_nI$ set classes are similarly
defined, but with the addition of flipping over; $T_n/T_nI$ set classes are
the union of pairwise $T_n$ and $T_nI$ set classes.   
The problem is: how many colour sets and
$T_n/T_nI$ set classes do exist?   For sake of clarity, let attention
be restricted to a simple case, $N=4$ and $L=3$: blue (B), green (G), red (R),
say e.g., \cite{Hoo07} \S47.

With regard to distinct colour sets and $T_n$ set classes, the results
coincide with their counterparts of the birthday-cake problem, where the
total number equals 81 and 24, respectively.   The total number of $T_nI$ set
classes, by definition, also equals 24.   What remains to be done is the
enumeration of repeating colour sets and $T_n/T_nI$ set classes.

Among the eight colour sets within a generic $T_n/T_nI$ set class, two
quadruplets exhibit colour sets related via circular permutation and, in
addition, four doublets related via reflection.   Accordingly, palindrome and
pseudo palindrome colour sets must be counted as repeating.    More
specifically,
$T_n$ and $T_nI$ i.e. pairwise $T_n$ set classes made of distinct colour sets
yield a single $T_n/T_nI$ set class, while single $T_n$ set classes made of
palindrome or pseudo palindrome colour sets yield a single $T_n/T_nI$ set
class.

For $x^4$ colour sets, $\{x,x,x,x\}$, the number of distinct 4-tuples is 
$4!/(4\cdot3\cdot2)=1$, and the number of different elements, $x$, equals
${3\choose1}=3$.   Accordingly, the total number of distinct and repeating
4-tuples is $3\cdot1=3$ and $3\cdot3=9$, respectively.   But colour sets are
palindrome in the case under discussion, which implies 4 additional 4-tuples
for a total of $3\cdot(3+4) = 21$ repeating colour sets. Then the total number
of distinct + repeating colour sets is 3+21 = 24.

For $x^3y$ colour sets, $\{x,x,x,y\}$, the number of distinct 4-tuples is 
$4!/(3\cdot2)=4$, and the number of different pairs, $xy$, equals
${4\choose2}=6$.   Accordingly, the total number of distinct and repeating
4-tuples is $6\cdot4=24$ and $0\cdot4=0$, respectively.   But colour sets are
pseudo
palindrome in the case under discussion, which implies 4 additional 4-tuples
for a total of $6\cdot(0+4)=24$ repeating colour sets.   Then the total number
of distinct + repeating colour sets is 24+24 = 48.

For $x^2y^2$ colour sets, $\{x,x,y,y\}$, the number of distinct 4-tuples is 
$4!/(2\cdot2)=6$, and the number of different pairs, $xy$, equals
${3\choose2}=3$.   Accordingly, the total number of distinct and repeating
4-tuples is $3\cdot6=18$ and $3\cdot(1+1)=6$, respectively, the last due to
$\{x,y,x,y\}$ and $\{y,x,y,x\}$.   But colour sets are palindrome or pseudo
palindrome in the case under discussion, which implies $4+4=8$
additional 4-tuples for a total of $3\cdot(2+8) = 30$ repeating colour sets.
Then the total number of distinct + repeating colour sets is 18+30 = 48.

For $x^2yz$ colour sets, $\{x,x,y,z\}$, the number of distinct 4-tuples is 
$4!/2=12$, and the number of different triads, $xyz$, equals
${3\choose1}=3$.   Accordingly, the total number of distinct and repeating
4-tuples is $3\cdot12=36$ and $0\cdot12=0$, respectively.   But colour sets,
$\{x,y,x,z\}$, are
pseudo palindrome in the case under discussion, which implies
4 additional 4-tuples for a total of $3\cdot(0+4)=12$ repeating colour sets.
Then the total number of distinct + repeating colour sets is $36+12=48$.

More specifically, in the case under discussion the total number of colour
sets belonging to an assigned kind can be determined along the following
steps.
  
\begin{trivlist}%
\item[\hspace\labelsep{\bf (a)}]
Start from an assigned $T_n/T_nI$ set class  of the kind considered.
\item[\hspace\labelsep{\bf (b)}]
Exchange $y$ with $z$ (if both present) to get a distinct (if any) $T_n/T_nI$
set class.
\item[\hspace\labelsep{\bf (c)}]
Exchange $x$ with $y$ (only once in case of multiplicity) to get a distinct
(if any) $T_n/T_nI$ set class.
\end{trivlist}%

For $x^4$ colour sets, the application of the method yields:
\begin{lefteqnarray*}
&& \{x,x,x,x\}, \{x,x,x,x\}, \{x,x,x,x\}, \{x,x,x,x\}~~, \\
&& \{x,x,x,x\}, \{x,x,x,x\}, \{x,x,x,x\}, \{x,x,x,x\}~~;
\end{lefteqnarray*}
where no additional exchange has to be performed as both $y$ and $z$ are
absent, yielding 1 $T_n/T_nI$ set class made of 1 distinct and 3+4 = 7
repeating colour
sets.   Keeping in mind $x={\rm B,G,R}$, there are $3\cdot1=3$ $T_n/T_nI$ set
classes for a total of $1\cdot3=3$ distinct and $7\cdot3=21$ repeating colour
sets, according to the above results.

For $x^3y$ colour sets, the application of the method yields:
\begin{lefteqnarray*}
&& \{x,x,x,y\}, \{x,x,y,x\}, \{x,y,x,x\}, \{y,x,x,x\}~~, \\
&& \{y,x,x,x\}, \{x,y,x,x\}, \{x,x,y,x\}, \{x,x,x,y\}~~;
\end{lefteqnarray*}
where no additional exchange has to be performed as $z$ is absent, yielding 1
$T_n/T_nI$ set class made of 4 distinct and 4 repeating colour
sets.   Keeping in mind $xy={\rm BG,GR,RB,GB,RG,BR}$, ($xy\ne yx$ due to 
different occurrences in $x$ and $y$), there are $6\cdot1=6$ $T_n/T_nI$ set
classes for a total of $6\cdot4=24$ distinct and $6\cdot4=24$ repeating colour
sets, according to the above results.

For $x^2y^2$ colour sets, the application of the method yields:
\begin{lefteqnarray*}
&& \{x,y,y,x\}, \{y,y,x,x\}, \{y,x,x,y\}, \{x,x,y,y\}~~, \\
&& \{x,y,y,x\}, \{x,x,y,y\}, \{y,x,x,y\}, \{y,y,x,x\}~~; \\
&& \{x,y,x,y\}, \{y,x,y,x\}, \{x,y,x,y\}, \{y,x,y,x\}~~, \\
&& \{y,x,y,x\}, \{x,y,x,y\}, \{y,x,y,x\}, \{x,y,x,y\}~~;
\end{lefteqnarray*}
where no additional exchange has to be performed as $z$ is absent, yielding 2
$T_n/T_nI$ set classes made of 4+2 = 6 distinct and 4+6 = 10 repeating colour
sets.   Keeping in mind $xy={\rm BG,GR,RB}$, ($xy=yx$ due to equal occurrences
in $x$ and $y$), there are $3\cdot2=6$ $T_n/T_nI$ set classes for a total of
$3\cdot(4+2)=18$ distinct and $3\cdot(4+6)=30$ repeating
colour sets, according to the above results.

For $x^2yz$ colour sets, the application of the method yields:
\begin{lefteqnarray*}
&& \{x,x,y,z\}, \{x,y,z,x\}, \{y,z,x,x\}, \{z,x,x,y\}~~, \\
&& \{z,y,x,x\}, \{x,z,y,x\}, \{x,x,z,y\}, \{y,x,x,z\}~~; \\
&& \{x,x,z,y\}, \{x,z,y,x\}, \{z,y,x,x\}, \{y,x,x,z\}~~, \\
&& \{y,z,x,x\}, \{x,y,z,x\}, \{x,x,y,z\}, \{z,x,x,y\}~~; \\
&& \{x,y,x,z\}, \{y,x,z,x\}, \{x,z,x,y\}, \{z,x,y,x\}~~, \\
&& \{z,x,y,x\}, \{x,z,x,y\}, \{y,x,z,x\}, \{x,y,x,z\}~~;
\end{lefteqnarray*}
where no additional exchange has to be performed, yielding 2 $T_n/T_nI$ set
classes (the first two coincide) made of $3\cdot4=12$ distinct and
$1\cdot4=4$ repeating colour
sets.   Keeping in mind $x={\rm B,G,R}$, $x\ne y\ne z$,
there are $3\cdot2=6$ $T_n/T_nI$ set
classes for a total of $3\cdot12=36$ distinct and $3\cdot4=12$ repeating
colour sets, according to the above results.

Finally, the total number of distinct, $N_{\rm C}$, repeating, $\Delta N$, and
total, $N_{\rm M}$, colour sets, respectively, reads:
\begin{lefteqnarray}
\label{eq:OCb}
&& N_{\rm C}=N_{\rm C}(x^4)+N_{\rm C}(x^3y)+N_{\rm C}(x^2y^2)+N_{\rm C}(x^2yz)
=3+24+18+36 \nonumber \\
&& \phantom{N_{\rm C}=N_{\rm C}(x^4)+N_{\rm C}(x^3y)+N_{\rm C}(x^2y^2)+
N_{\rm C}(x^2yz)}=81~~; \\
\label{eq:DOb}
&& \Delta N=\Delta N(x^4)+\Delta N(x^3y)+\Delta N(x^2y^2)+\Delta N(x^2yz)=
21+24+30+12 \nonumber \\
&& \phantom{\Delta N=\Delta N(x^4)+\Delta N(x^3y)+\Delta N(x^2y^2)+\Delta N
(x^2yz)}=87~~; \\
\label{eq:OMb}
&& N_{\rm M}=N_{\rm C}+\Delta N=24+48+48+48=168~~; 
\end{lefteqnarray}
and the number of $T_n/T_nI$ set classes, $\nu_{\rm M}$, is:
\begin{lefteqnarray}
\label{eq:oMb}
&& \nu_{\rm M}=\frac{N_{\rm M}}8=\frac{24}8+\frac{48}8+\frac{48}8+\frac{48}8=
3+6+6+6=21~~;
\end{lefteqnarray}
or 168/8=21.

Within the framework of group theory, the application of Burnside's lemma to
necklace problem yields e.g., \cite{Hoo07} \S47:
\begin{lefteqnarray}
\label{eq:oMB}
&& \nu_{\rm M}=\frac{84}8+\frac{12}8+\frac{36}8+\frac{36}8=\frac{168}8=21~~;
\end{lefteqnarray}
which is in agreement with Eq.\,(\ref{eq:oMb}) but the addends are different.

Within the framework of group theory, the application of Polya's theorem to
necklace problem yields e.g., \cite{Hoo07} \S64:
\begin{lefteqnarray}
\label{eq:oMP}
&& \nu_{\rm M}=(1+1+1)+(1+1+1+1+1+1)+2(1+1+1)+2(1+1+1) \nonumber \\
&& \phantom{\nu_{\rm M}}=3+6+6+6=21~~;
\end{lefteqnarray}
which is in agreement with Eq.\,(\ref{eq:oMb}), exhibiting equal addends.

The result from Eq.\,(\ref{eq:oMB}) is based on Burnside's lemma and does not
involve weights on $T_n/T_nI$ set classes.   On the other hand, the weighted
method given by Polya's theorem is simply a refinement of the same idea,
showing exactly how the 21 $T_n/T_nI$ set classes are distributed among the
various weights e.g., \cite{Hoo07} \S63.   For further details, an interested reader is
addressed to specialized monographs e.g., \cite{Hoo07} or textbooks e.g.,
\cite{Ben06} Chap.\,9.   The current method exactly reproduces the results
inferred via Polya's theorem.

\section{Basic ideas on $\Psi_{12}^n$ $n$-hedrons}
\label{a:nhed}
\subsection{Analytic geometry in $\Re^k$}
\label{a:angk}

The formulation of analytic geometry in ordinary Euclidean space, $\Re^3$, can
easily be extended to Euclidean $k$-spaces, $\Re^k$.

With regard to a Cartesian orthogonal reference frame,
$({\sf O}\,x_1\,x_2,...\,x_k)$, let 
${\sf P_A}\equiv(x_{A1},x_{A2},...,x_{Ak})$ and
${\sf P_B}\equiv(x_{B1},x_{B2},...,x_{Bk})$ be generic points.   The related
square distance reads:
\begin{equation}
\label{eq:PAB}
({\sf \overline{P_AP_B}})^2=\sum_{i=1}^k(x_{Ai}-x_{Bi})^2~~;
\end{equation}
if ${\sf P_A}$ is fixed (the origin say) while ${\sf P_B}$ remains generic,
then Eq.\,(\ref{eq:PAB}) represents a $k$-dimension sphere, or $k$-sphere, of
radius, $R={\sf \overline{P_AP_B}}$, centered on ${\sf P_A}$.

The mean point, ${\sf P_M}$, of the segment joining the points, ${\sf P_A}$,
${\sf P_B}$, reads:
\begin{equation}
\label{eq:PM}
{\sf P_M}\equiv\left(\frac{x_{A1}+x_{B1}}2,\frac{x_{A2}+x_{B2}}2,...,
\frac{x_{Ak}+x_{Bk}}2\right)~~;
\end{equation}
and the identities, ${\sf \overline{P_AP_M}}={\sf \overline{P_BP_M}}$,
${\sf \overline{P_AP_M}}+{\sf \overline{P_BP_M}}={\sf \overline{P_AP_B}}$, can
easily be verified via Eqs.\,(\ref{eq:PAB}) and (\ref{eq:PM}).

Let $r$, $r^\prime$, be generic straight lines in $\Re^k$, expressed as:
\begin{lefteqnarray}
\label{eq:r}
&& r:\qquad
\frac{x_1-x_{10}}{L_1}=\frac{x_2-x_{20}}{L_2}=...=\frac{x_k-x_{k0}}{L_k}~~;
\\
\label{Li}
&& \phantom{r:\qquad}L_i=\overline x_{i0}-x_{i0}~~;\qquad 1\le i\le k~~; \\
\label{eq:rp}
&& r^\prime:\qquad
\frac{x_1^\prime-x_{10}^\prime}{L_1^\prime}=\frac{x_2^\prime-x_{20}^\prime}
{L_2^\prime}=...=\frac{x_k^\prime-x_{k0}^\prime}{L_k^\prime}~~; \\
\label{Lip}
&& \phantom{r^\prime:\qquad}
L_i^\prime=\overline x_{i0}^\prime-x_{i0}^\prime~~;\qquad 1\le i\le k~~;
\end{lefteqnarray}
where ${\sf P_0}\equiv(x_{10},x_{20},...,x_{k0})$, ${\sf\overline P_0}\equiv
(\overline x_{10},\overline x_{20},...,\overline x_{k0})$, are selected points
on $r$ and  ${\sf P_0^\prime}\equiv(x_{10}^\prime,x_{20}^\prime,...,
x_{k0}^\prime)$, ${\sf\overline P_0^\prime}\equiv(\overline x_{10}^\prime,
\overline x_{20}^\prime,...,\overline x_{k0}^\prime)$, are selected points
on $r^\prime$. 

The angle, $\widehat{rr^\prime}$, can be expressed as:
\begin{equation}
\label{eq:rrp}
\vert\cos\widehat{rr^\prime}\vert=\left\vert\displayfrac{\sum_{i=1}^kL_iL_i^
\prime}
{\left[\sum_{i=1}^kL_i^2\right]^{1/2}\left[\sum_{i=1}^k(L_i^\prime)^2\right]^
{1/2}}\right\vert~~;
\end{equation}
where the orthogonality condition reads $\cos\widehat{rr^\prime}=0$.

Let $p$ be a generic $(k-1)$-plane in $\Re^k$, expressed as:
\begin{lefteqnarray}
\label{eq:p}
&& p:\qquad a_1x_1+a_2x_2+...+a_kx_k+a_0=0~~;
\end{lefteqnarray}
where $-a_0/a_i$, is the intercept of $p$ on the coordinate axis, $x_i$,
$1\le i\le k$. 

The angle, $\widehat{rp}$, between the straight line, $r$, and the
$(k-1)$-plane, $p$, can be expressed as:
\begin{equation}
\label{eq:arp}
\vert\sin\widehat{rp}\vert=\left\vert\displayfrac{\sum_{i=1}^kL_ia_i}
{\left[\sum_{i=1}^kL_i^2\right]^{1/2}\left[\sum_{i=1}^ka_i^2\right]^{1/2}}
\right\vert~~;
\end{equation}
where the orthogonality condition reads $\sin\widehat{rp}=\mp1$.

For further details, an interested reader is addressed to earlier
investigations \cite{Cai00} Chap.\,2 \S2.16 \cite{Cai13} \cite{Cai15} \cite{Cai16}
Chap.\,4 \S4.17.

\subsection{General properties of $(k+1)$-hedrons}
\label{a:gpkh}

With regard to an Euclidean $k$-space, $\Re^k$, and a Cartesian orthogonal
reference frame, $({\sf O}x_1\,x_2\,...\,x_k)$, let $(k+1)$ $k$-misaligned
points be assigned as ${\sf V_1}$, ${\sf V_2}$, ..., ${\sf V_{k+1}}$.   Points
are $k$-misaligned in the sense different $k$-tuples,
$\{{\sf V_{i_1}}, {\sf V_{i_2}}, ..., {\sf V_{i_k}}\}$,
$1\le i_1<i_2<...<i_k\le k+1$, belong to different $(k-1)$-planes in $\Re^k$.
Let each point be connected with the remaining ones by segments.   Let
the resulting geometrical figure be defined as $(k+1)$-hedron, and denoted as
$\Phi^{k+1}$.   In ordinary space, $k=3$, 4-hedrons reduce to ordinary
tetrahedrons.

The above definition could appear ambiguous, in that the etymology of
tetrahedron implies the existence of four faces.   Keeping in mind
$(k+1)$-hedrons have $k$-dimensions, it can be seen (as shown below) that the
number of $(k-1)$-faces equals the number of connected points, $(k+1)$.   Then
etymological meaning is preserved provided $(k-1)$-faces are considered
instead of ordinary (2-) faces.

Another source of ambiguity lies in the fact, that ``$(k+1)$-hedrons'' are
already defined in ordinary space, $\Re^3$, concerning ordinary faces e.g.,
hexahedron, octahedron, dodecahedron, icosahedron.   Throughout the present
paper, according to the current definition, $(k+1)$-hedrons are intended as
solids in $\Re^k$ where a generic point,  ${\sf V_i}$, $1\le i\le k+1$, is
connected to the remaining ($k$-misaligned) $k$.

Any point, ${\sf V_i}$, $1\le i\le k+1$, is a vertex or 0-face of
$\Phi^{k+1}$.   Accordingly, the number of $0$-faces equals
the number of distinct $1$-tuples of vertexes, ${k+1 \choose 1}=(k+1)/1!$.

Any duo of vertexes, ${\sf V_iV_j}$, $i\ne j$, yields a side or 1-face of
$\Phi^{k+1}$.   Accordingly, the number of 1-faces equals the number of
distinct 2-tuples of vertexes, ${k+1 \choose 2}=(k+1)k/2!$.

Any trio of vertexes, ${\sf V_iV_jV_\ell}$, $i\ne j\ne\ell$, yields a triangle
or 2-face of $\Phi^{k+1}$.   Accordingly, the number of 2-faces equals the
number of distinct 3-tuples of vertexes, ${k+1 \choose 3}=(k+1)k(k-1)/3!$.

Any quartet of vertexes, ${\sf V_iV_jV_\ell  V_m}$, $i\ne j\ne\ell\ne m$,
yields a tetrahedron or 3-face of $\Phi^{k+1}$.   Accordingly, the number of
3-faces equals the number of distinct 4-tuples of vertexes,
${k+1 \choose 4}=(k+1)k(k-1)(k-2)/4!$.

In general, any $j$-tet of vertexes, ${\sf V_{i_1}V_{i_2}...V_{i_j}}$,
$i_1\ne i_2\ne...\ne i_j$, $1\le j\le k$, yields a $j$-hedron or
$(j-1)$-face of $\Phi^{k+1}$.    Accordingly, the number of $(j-1)$-faces
equals the number of distinct $j$-tuples of vertexes,
${k+1 \choose j}=(k+1)k...(k-j+2)/j!$.

To this respect, a useful recursion formula follows from the identity:
\begin{equation}
\label{eq:focb}
{k+1 \choose j}={k \choose j}+{k \choose j-1}~~;\qquad1\le j\le k~~;
\end{equation}
which translates into:
\begin{equation}
\label{eq:fofa}
F_j(\Phi^{k+1})=F_j(\Phi^k)+F_{j-1}(\Phi^k)~~;\qquad1\le j\le k~~;
\end{equation}
where, in general, $F_\ell(\Phi^{m+1})$ is the number of $(\ell-1)$-faces of
$\Phi^{m+1}$, $1\le\ell\le m$.

\subsection{A special case: $\Psi_{12}^n$ regular inclined $n$-hedrons}
\label{a:psi12n}

With regard to an Euclidean $n$-space, $\Re^n$, and a Cartesian orthogonal
reference frame, $({\sf O}x_1\,x_2\,...\,x_n)$, let a regular, inclined,
$n$-hedron where vertexes are placed on the coordinate axes at a
distance, $x_i=12$, $1\le i\le12$, from the origin, be considered and
denoted as $\Psi_{12}^n$.   Points with positive integer coordinates, or
positive integer points, within $\Psi_{12}^n$ satisfy Eq.\,(\ref{eq:boco}) as
shown in the text.   More specifically, the extension of Eq.\,(\ref{eq:boco})
to real coordinates reads:
\begin{lefteqnarray}
\label{eq:pn}
&& p_n:\qquad x_1+x_2+...+x_n=12~~;\qquad1\le n\le12~~;
\end{lefteqnarray}
which represents a $(n-1)$-plane where the intercepts on the coordinate axes
coincide with the vertexes of $\Psi_{12}^n$.

The $n$-sector ($n=2$, bisector; $n=3$, trisector; and so on) of the positive
$2^n$-ant, by definition, includes both the origin and the point of unit
coordinates, which implies $x_{i0}=0$, $\overline x_{i0}=1$, $1\le i\le n$,
and Eq.\,(\ref{eq:r}) reduces to:
\begin{lefteqnarray}
\label{eq:rn}
&& r_n:\qquad x_1=x_2=...=x_n~~;\qquad1\le n\le12~~;
\end{lefteqnarray}
which is the locus of points with identical coordinates.

The angle between the $(n-1)$-plane, $p_n$, and the $n$-sector, $r_n$, keeping
in mind $a_i=1$, $L_i=1$, $1\le i\le n$, via Eq.\,(\ref{eq:arp}) reads:
\begin{equation}
\label{eq:rnpn}
\vert\sin\widehat{r_np_n}\vert=\left\vert\frac n{\sqrt n\sqrt n}\right\vert=1
~~;
\end{equation}
which implies $r_n$ is orthogonal to $p_n$ through the orthocentre of
$\Psi_{12}^n$ that is, in fact, the sole point of $p_n$ with identical
coordinates.

With regard to the generic $(j-1)$-face, $1\le j\le n$, belonging to
$\Psi_{12}^n$, the
orthocentre, ${\sf H_j}$, exhibits $j$ equal coordinates in that it is
equally distant from $j$ vertexes of $\Psi_{12}^n$.   Accordingly, nonzero
coordinates equal $12/j$ via Eq.\,(\ref{eq:boco}).   The result is:
\begin{lefteqnarray}
\label{eq:Hj1}
&& {\sf H_j}(i_1,i_2,...,i_j)\equiv\left[\frac{12}j\sum_{k=1}^j
\delta_{i_k1},\frac{12}j\sum_{k=1}^j\delta_{i_k2},...,\frac{12}j\sum_{k=1}^j
\delta_{i_kj}\right]~~; \nonumber \\
&& \phantom{{\sf H_{j-1}}(i_1,i_2,...,i_j)\equiv~}1\le i_1<i_2<...<i_j\le n~~;
\end{lefteqnarray}
where regular vertexes, regular sides, regular triangles, regular
tetrahedrons, and so
on, are conceived as 0-faces, 1-faces, 2-faces,  3-faces, and so on,
respectively.   In particular, the orthocentre of $\Psi_{12}^n$, conceived as
related to a single $(n-1)$-face, reads:
\begin{lefteqnarray}
\label{eq:Hn1}
&& {\sf H_n}(1,2,...,n)\equiv\left(\frac{12}n,\frac{12}n,...,\frac{12}n
\right)~~;
\end{lefteqnarray}
where ${\sf H_n}$ is a positive integer point only if the ratio, 12/$n$,
is integer, which implies $n=1,2,3,4,6,12$.

The number of $(j-1)$-faces, $F_j$, $1\le j\le n$, belonging to
$\Psi_{12}^n$, $1\le n\le13$, can be inferred from the general results
mentioned above, where $k=n-1$, as listed in Table \ref{t:face}.   More
specifically, defining $(-1)$-faces and $(n-1)$-faces as $n$-hedron
metacentres and $n$-hedrons, respectively, related values equal unity and
Table \ref{t:face} reproduces Pascal's triangle e.g., \cite{Hoo07} \S30.

The metacentre of $\Psi_{12}^n$ is defined as the centre of a
$n$-sphere, on the surface of which vertexes of $\Psi_{12}^n$ are
lying i.e. the origin of coordinates.   In addition, the metacentre can be
thought of as the orthocentre of a $(n+1)$-hedron where $\Psi_{12}^n$ is a
$(n-1)$-face and the additional vertex lies on the $n$-sector of the
negative $2^n$-ant, at a distance, $R_n=12$, from the origin.

 In the special case, $n=0$,
related 0-hedron lies outside  $\Re^0$ and no vertex can be defined, while the
metacentre coincides with $\Re^0$, yielding 1 $(-1)$-face.
\begin{table}
\caption{Number of $(j-1)$-faces, $F_j$, $0\le j\le n$, belonging to
$\Psi_{12}^n$, $0\le n\le13$.   In particular, $n=F_1$.   The addition of
$(-1)$-faces, $F_0$, related to $n$-hedron metacentres, and $(n-1)$-faces,
$F_n$, related to $n$-hedrons, makes the horizontal and oblique unit line,
respectively, yielding Pascal's triangle.   For $j<10$, $F_j$ is denoted as
$F_{0j}$ to save aesthetics.   See text for further details.}
\label{t:face}
\begin{center}
\begin{tabular}{|r|r|r|r|r|r|r|r|r|r|r|r|r|r|r|} \hline
\hline
$F_{00}$ & 1 & 1         &  1 &  1                  &  1                  &   1                  &   1 &   1 &   1 &   1                  &   1                 &   1                 &   1 &    1                     \\
$F_{01}$ &   & 1         &  2 &  3                  &  4                  &   5                  &   6 &   7 &   8 &   9                  &  10                 &  11                 &  12 &   13                     \\
$F_{02}$ &   &           &  1 &  3                  &  6                  &  10                  &  15 &  21 &  28 &  36                  &  45                 &  55                 &  66 &   78                     \\
$F_{03}$ &   &           &    &  1                  &  4                  &  10                  &  20 &  35 &  56 &  84                  & 120                 & 165                 & 220 &  286                     \\
$F_{04}$ &   &           &    &                     &  1                  &   5                  &  15 &  35 &  70 & 126                  & 210                 & 330                 & 495 &  715                     \\
$F_{05}$ &   &           &    &                     &                     &   1                  &   6 &  21 &  56 & 126                  & 252                 & 462                 & 792 & 1287                     \\
$F_{06}$ &   &           &    &                     &                     &                      &   1 &   7 &  28 &  84                  & 210                 & 462                 & 924 & 1716                     \\
$F_{07}$ &   &           &    &                     &                     &                      &     &   1 &   8 &  36                  & 120                 & 330                 & 792 & 1716                     \\
$F_{08}$ &   &           &    &                     &                     &                      &     &     &   1 &   9                  &  45                 & 165                 & 495 & 1287                     \\
$F_{09}$ &   &           &    &                     &                     &                      &     &     &     &   1                  &  10                 &  55                 & 220 &  715                     \\
$F_{10}$ &   &           &    &                     &                     &                      &     &     &     &                      &   1                 &  11                 &  66 &  286                     \\
$F_{11}$ &   &           &    &                     &                     &                      &     &     &     &                      &                     &   1                 &  12 &   78                     \\
$F_{12}$ &   &           &    &                     &                     &                      &     &     &     &                      &                     &                     &   1 &   13                     \\
$F_{13}$ &   &           &    &                     &                     &                      &     &     &     &                      &                     &                     &     &    1                     \\
\hline                            
\end{tabular}                     
\end{center}                      
\end{table}                       

Let a $(j-1)$-face and a $(n-j-1)$-face of $\Psi_{12}^n$ be assigned.   The
coordinates of generic points lying within, ${\sf P_j}$,
${\sf P_{n-j}}$, respectively, read:
\begin{leftsubeqnarray}
\slabel{eq:Pjna}
&& {\sf P_j}\equiv\left(\ell_1\sum_{k=1}^j
\delta_{i_k1},\ell_2\sum_{k=1}^j\delta_{i_k2},...,\ell_n\sum_{k=1}^j
\delta_{i_kj}\right)~~; \nonumber \\
&& \phantom{{\sf H_{j-1}}\equiv~}1\le i_1<i_2<...<i_j\le n~~; \\
\slabel{eq:Pjnb}
&& {\sf P_{n-j}}\equiv\left(\ell_1^\prime\sum_{m=1}^{n-j}
\delta_{i_m^\prime1},\ell_2^\prime\sum_{m=1}^{n-j}\delta_{i_m^\prime2},...,
\ell_n^\prime\sum_{m=1}^{n-j}\delta_{i_m^\prime j}\right)~~; \nonumber \\
&& \phantom{{\sf H_{n-1-j}}\equiv~}1\le i_1^\prime<i_2^\prime<...<i_{n-j}^
\prime\le n~~;
\label{seq:Pjn}
\end{leftsubeqnarray}
and the hyperfaces under discussion are defined as opposite if, in addition,
$i_k\ne i_m^\prime$, $1\le k\le j$, $1\le m\le n-j$.

The special case, ${\sf P_j}\equiv{\sf H_j}$,
${\sf P_{n-j}}\equiv{\sf H_{n-j}}$, via Eq.\,(\ref{eq:Hj1}) reads:
\begin{leftsubeqnarray}
\slabel{eq:Hjna}
&& {\sf H_j}\equiv\left(\frac{12}j\sum_{k=1}^j
\delta_{i_k1},\frac{12}j\sum_{k=1}^j\delta_{i_k2},...,\frac{12}j\sum_{k=1}^j
\delta_{i_kn}\right)~~; \nonumber \\
&& \phantom{{\sf H_{j-1}}\equiv~}1\le i_1<i_2<...<i_j\le n~~; \\
\slabel{eq:Hjnb}
&& {\sf H_{n-j}}\equiv\left(\frac{12}{n-j}\sum_{m=1}^{n-j}\delta_
{i_m^\prime1},\frac{12}{n-j}\sum_{m=1}^{n-j}\delta_{i_m^\prime2},...,\frac{12}
{n-j}\sum_{m=1}^{n-j}\delta_{i_m^\prime n}\right)~~; \nonumber \\
&& \phantom{{\sf H_{n-1-j}}\equiv~}1\le i_1^\prime<i_2^\prime<...<i_{n-j}^
\prime\le n~~;
\label{seq:Hjn}
\end{leftsubeqnarray}
where $i_k\ne i_m^\prime$; $1\le k\le j$; $1\le m\le n-j$; and the orthocentre
of $\psi_{12}^n$, ${\sf H_n}$, is necessarily aligned with ${\sf H_j}$
and ${\sf H_{n-j}}$.

More specifically, the particularization of Eq.\,(\ref{eq:r}) to
$x_{i0}=(12/j)(\delta_{i_1i}+\delta_{i_2i}+...+\delta_{i_ji})$,
$\overline x_{i0}=[12/(n-j)](\delta_{i_1^\prime i}+\delta_{i_2^\prime i}+...+
\delta_{i_{n-j}^\prime i})$, $1\le i\le n$, yields:
\begin{lefteqnarray}
\label{eq:HHjn}
&& \displayfrac{x_1-\frac{12}j\sum_{k=1}^j\delta_{i_k1}}{\frac{12}{n-j}
\sum_{m=1}^{n-j}\delta_{i_m^\prime1}-\frac{12}j\sum_{k=1}^j\delta_{i_k1}}=
\displayfrac{x_2-\frac{12}j\sum_{k=1}^j\delta_{i_k1}}{\frac{12}{n-j}
\sum_{m=1}^{n-j}\delta_{i_m^\prime1}-\frac{12}j\sum_{k=1}^j\delta_{i_k1}}=...
\nonumber \\
&& \phantom{\displayfrac{x_1-\frac{12}j\sum_{k=1}^j\delta_{i_k1}}{\frac{12}
{n-j}\sum_{m=1}^{n-j}\delta_{i_m^\prime1}-\frac{12}j\sum_{k=1}^j\delta_{i_k1}}
}
=\displayfrac{x_n-\frac{12}j\sum_{k=1}^j\delta_{i_kn}}{\frac{12}{n-j}
\sum_{m=1}^{n-j}\delta_{i_m^\prime n}-\frac{12}j\sum_{k=1}^j\delta_{i_kn}}~~;
\end{lefteqnarray}
where, in the case under discussion of opposite hyperfaces, each sum equals
zero or
unity and attains different values on different points.   Accordingly, the
generic term appearing in Eq.\,(\ref{eq:HHjn}) reads either $-jx_i/12+1$, for
a total of $j$ occurrences, or $(n-j)x_i/12$, for a total of $(n-j)$
occurrences, $1\le i\le n$.

In the special case of the orthocentre of $\Psi_{12}^n$, $x_i=12/n$,
$1\le i\le n$, and the above mentioned alternatives reduce to either
$-(j/12)(12/n)+1=(n-j)/n$ or $[(n-j)/12](12/n)=(n-j)n$, which implies
Eq.\,(\ref{eq:HHjn}) is satisfied and, in turn, the orthocentre of
$\Psi_{12}^n$ is aligned with the orthocentre of opposite $(j-1)$ and
$(n-j-1)$-faces.

Let ${\sf P_n}\equiv(\ell_1,\ell_2,...,\ell_n)$ be a generic integer point
belonging to $\Psi_{12}^n$.   The square distance,
$(\overline{\sf P_nH_n})^2=R_n^2$, from the orthocentre of
$\Psi_{12}^n$, via Eqs.\,(\ref{eq:PAB}) and (\ref{eq:Hn1}) reads:
\begin{equation}
\label{eq:Rn1}
R_n^2=\sum_{i=1}^n(\ell_i-\ell_0)^2~~;\qquad\ell_0=\frac{12}n~~;
\end{equation}
where $R_n$ is the radius of a $(n-1)$-circle centered on ${\sf H_n}$.

The straight line joining ${\sf P_n}$ and ${\sf H_n}$, via
Eq.\,(\ref{eq:r}) is expressed as:
\begin{equation}
\label{eq:ra}
r_a:\qquad\frac{x_1-\ell_0}{\ell_1-\ell_0}=\frac{x_2-\ell_0}{\ell_2-\ell_0}=
...=\frac{x_n-\ell_0}{\ell_n-\ell_0}~~;
\end{equation}
where $L_k=\ell_k-\ell_0$, $1\le k\le n$.

The counterpart of Eq.\,(\ref{eq:ra}), related to a selected permutation of
the coordinates of ${\sf P_n}$, is:
\begin{equation}
\label{eq:rb}
r_b:\qquad\frac{x_1-\ell_0}{\ell_{i_1}-\ell_0}=\frac{x_2-\ell_0}
{\ell_{i_2}-\ell_0}=...=\frac{x_n-\ell_0}{\ell_{i_n}-\ell_0}~~;
\end{equation}
where $L_{i_k}=\ell_{i_k}-\ell_0$, $1\le k\le n$.

The angle, $\widehat{r_ar_b}$, via Eq.\,(\ref{eq:rrp}) reads:
\begin{equation}
\label{eq:rab}
\vert\cos\widehat{r_ar_b}\vert=\left\vert\displayfrac{\sum_{k=1}^n
(\ell_k-\ell_0)(\ell_{i_k}-\ell_0)}
{\left[\sum_{i=1}^k(\ell_k-\ell_0)^2\right]^{1/2}
\left[\sum_{i=1}^k(\ell_{i_k}-\ell_0)^2\right]^{1/2}}\right\vert~~;
\end{equation}
where, in the case under discussion, the sums of squares exhibit same
addends in different order.   Accordingly, Eq.\,(\ref{eq:rab}) reduces to:
\begin{equation}
\label{eq:abp}
\vert\cos\widehat{r_ar_b}\vert=\left\vert\displayfrac{\sum_{k=1}^n
(\ell_k-\ell_0)(\ell_{i_k}-\ell_0)}
{\sum_{i=1}^k(\ell_k-\ell_0)^2}\right\vert~~;
\end{equation}
where the permutation of the coordinates acts only on the sum of products.

The special case of circular permutation, $i_k=k+1$, reads:
\begin{equation}
\label{eq:abc}
\vert\cos\widehat{r_ar_b}\vert=\left\vert\displayfrac{\sum_{k=1}^n
(\ell_k-\ell_0)(\ell_{k+1}-\ell_0)}
{\sum_{i=1}^k(\ell_k-\ell_0)^2}\right\vert~~;
\end{equation}
where $\ell_{n+1}=\ell_1$.   It is apparent adjacent circular permutations of
coordinates e.g., $(\ell_1,\ell_2,...,\ell_n)$, $(\ell_2,\ell_3,...,\ell_1)$;
$(\ell_2,\ell_3,...,\ell_1)$, $(\ell_3,\ell_4,...,\ell_2)$; yield the same
angle in that the sum of products in Eq.\,(\ref{eq:abc}) exhibits same
addends in different order.    On the other hand, non adjacent circular
permutations of coordinates yield different angles with respect to the sum of
adjacent angles, as related straight lines
are not complanar.   In fact, positive integer points whose coordinates belong
to a same $T_n$ set class are similarly placed within different
$\Psi_{12,i}^n$, $1\le i\le n$, i.e. congruent $(n-1)$-hedrons, where the
$i$th vertex coincides with the orthocentre of $\Psi_{12}^n$ and the remaining
ones coincide with vertexes (different from the $i$th) of $\Psi_{12}^n$.

The special case of reflection, $i_k=n-k+1$, reads:
\begin{equation}
\label{eq:abr}
\vert\cos\widehat{r_ar_b}\vert=\left\vert\displayfrac{\sum_{k=1}^n
(\ell_k-\ell_0)(\ell_{n-k+1}-\ell_0)}
{\sum_{i=1}^k(\ell_k-\ell_0)^2}\right\vert~~;
\end{equation}
where adjacent circular permutations yield the same angle as the sum of
products in Eq.\,(\ref{eq:abr}) exhibits same addends in different order.
In fact, positive integer points with coordinates belonging to a same
$T_n/T_nI$ set class are similarly and symmetrically placed within different
$\Psi_{12,i}^n$, $1\le i\le n$.

The mean point, ${\sf M_n}$, of the segment joining pairwise points,
${\sf P_n}\equiv(\ell_1,\ell_2,...,\ell_n)$,
${\sf Q_n}\equiv(\ell_n,\ell_{n-1},...,\ell_1)$, via Eq.\,(\ref{eq:PM}) reads:
\begin{equation}
\label{eq:Mn}
{\sf M_n}\equiv\left(\frac{\ell_1+\ell_n}2,\frac{\ell_2+\ell_{n-1}}2,...,
\frac{\ell_n+\ell_1}2,\right)~~;
\end{equation}
accordingly, the coordinates of ${\sf M_n}$ are palindrome as either
$(a_1,...,a_{n/2},a_{n/2},...,a_1)$ or
$(a_1,...,a_{(n-1)/2},a_{(n+1)/2},a_{(n-1)/2},...,a_1)$ according if $n$ is
even or odd, respectively.

To get further insight, let $\Psi_{12}^n$ be conceived as a series of
Chinese-box $n$-hedronical shells, $\psi_{k,12}^n$, where the external
shell, $\psi_{0,12}^n$, has not to be considered in that related points
exhibit one null coordinate at least.   The last shell of the series relates
to the value of $k$ for which the condition, $k\le12/n$, still holds.   The
special case, $k=12/n$, implies $\psi_{12/n,12}^n$ coincides with the
orthocentre of $\Psi_{12}^n$, provided $n$ is a divisor of 12.   For instance,
vertexes of $\psi_{k,12}^n$ read ${\sf V_i^\prime}\equiv[k+\delta_{i1}(12-kn),
k+\delta_{i2}(12-kn),...,k+\delta_{in}(12-kn)]$, and mean points between
vertexes read ${\sf M_{ij}^\prime}\equiv[k+(\delta_{i1}+\delta_{j1})(12-kn)/2,
k+(\delta_{i2}+\delta_{j2})(12-kn)/2,...,k+(\delta_{in}+\delta_{jn})
(12-kn)/2]$.

For even $n$, the straight line joining mean points on opposite side of
$\psi_{k,12}^n$, restricting to palindrome coordinates, via Eq.\,(\ref{eq:r})
is expressed as:
\begin{lefteqnarray}
\label{eq:H2H2}
&& \frac{x_{n/2-1}-k}{[k+(12-kn)/2]-k}=\frac{x_{n/2}-[k+(12-kn)/2]}
{k-[k+(12-kn)/2]} \nonumber \\
&& =\frac{x_{n/2+1}-[k+(12-kn)/2]}{k-[k+(12-kn)/2]}=
\frac{x_{n/2+2}-k}{[k+(12-kn)/2]-k}~~;
\end{lefteqnarray}
which, after some algebra, reduces to:
\begin{lefteqnarray}
\label{eq:H2H21}
&& x_{n/2-1}=-x_{n/2}+2k+6-\frac{kn}2=-x_{n/2+1}+2k+6-\frac{kn}2=x_{n/2+2}~~;
\qquad
\end{lefteqnarray}
while the remaining coordinates are fixed as:
\begin{lefteqnarray}
\label{eq:H2H2f}
&& x_1=...=x_{n/2-2}=x_{n/2+3}=...=x_n=k~~;
\end{lefteqnarray}
in connection with the selected $\psi_{k,12}^n$.

In conclusion, positive integer points exhibiting palindrome coordinates, for
even $n$ lie on the straight line, expressed by
Eqs.\,(\ref{eq:H2H21})-(\ref{eq:H2H2f}).   For instance, let $n=6$, $k=1$, be
considered.   Accordingly, Eqs.\,(\ref{eq:H2H21})-(\ref{eq:H2H2f}) reduce to:
\begin{lefteqnarray}
\label{eq:n6k1}
&& x_2=-x_3+5=-x_4+5=x_5~~;\qquad x_1=x_6=1~~;
\end{lefteqnarray}
where $(1,4,1,1,4,1)$, $(1,1,4,4,1,1)$, are coordinates of mean points on
opposite sides of related 4-face, and $(1,3,2,2,3,1)$, $(1,2,3,3,2,1)$, are
palindrome coordinates of positive integer points on the straight line joining
the above mentioned mean points, expressed by Eq.\,(\ref{eq:n6k1}).

For odd $n$, the straight line joining a vertex with the mean point of the
opposite sides of a 3-face of $\psi_{k,12}^n$, restricting to palindrome
coordinates, via Eq.\,(\ref{eq:r}) is expressed as:
\begin{lefteqnarray}
\label{eq:H0H3}
&& \frac{x_{(n-1)/2}-k}{[k+(12-kn)/2]-k}=\frac{x_{(n+1)/2}-[k+(12-kn)]}
{k-[k+(12-kn)/2]} \nonumber \\
&& \phantom{\frac{x_{(n-1)/2}-k}{[k+(12-kn)/2]-k}}
=\frac{x_{(n+3)/2}-k}{[k+(12-kn)/2]-k}~~;
\end{lefteqnarray}
which, after some algebra, reduces to:
\begin{lefteqnarray}
\label{eq:H0H31}
&& 2x_{(n-1)/2}=-x_{(n+1)/2}+12-kn+3k=2x_{(n+3)/2}~~;
\qquad
\end{lefteqnarray}
while the remaining coordinates are fixed as:
\begin{lefteqnarray}
\label{eq:H0H3f}
&& x_1=...=x_{(n-3)/2}=x_{(n+5)/2}=...=x_n=k~~;
\end{lefteqnarray}
in connection with the selected $\psi_{k,12}^n$.

In conclusion, positive integer points exhibiting palindrome coordinates, for
odd $n$ lie on the straight line, expressed by
Eqs.\,(\ref{eq:H0H31})-(\ref{eq:H0H3f}).   For instance, let $n=5$, $k=1$, be
considered.   Accordingly, Eqs.\,(\ref{eq:H0H31})-(\ref{eq:H0H3f}) reduce to:
\begin{lefteqnarray}
\label{eq:n5k1}
&& 2x_2=-x_3+10=2x_4~~;\qquad x_1=x_5=1~~;
\end{lefteqnarray}
where $(1,1,8,1,1)$, $(1,9/2,1,9/2,1)$, are coordinates of the vertex and mean
point on
opposite side of related 3-face, respectively, and $(1,2,6,2,1)$,
$(1,3,4,3,1)$, $(1,4,2,4,1)$, are
palindrome coordinates of positive integer points on the straight line joining
the above mentioned vertex and mean point, expressed by Eq.\,(\ref{eq:n5k1}).

\end{document}